\documentclass[11pt]{article}
\usepackage{amssymb,amsmath,latexsym}
\begin{document}

\begin{doublespace}

\newtheorem{thm}{Theorem}[section]
\newtheorem{lemma}[thm]{Lemma}
\newtheorem{defn}[thm]{Definition}
\newtheorem{prop}[thm]{Proposition}
\newtheorem{corollary}[thm]{Corollary}
\newtheorem{remark}[thm]{Remark}
\newtheorem{example}[thm]{Example}
\numberwithin{equation}{section}
\def\ee{\varepsilon}
\def\qed{{\hfill $\Box$ \bigskip}}
\def\NN{{\cal N}}
\def\AA{{\cal A}}
\def\MM{{\cal M}}
\def\BB{{\cal B}}
\def\CC{{\cal C}}
\def\LL{{\cal L}}
\def\DD{{\cal D}}
\def\FF{{\cal F}}
\def\EE{{\cal E}}
\def\QQ{{\cal Q}}
\def\RR{{\mathbb R}}
\def\R{{\mathbb R}}
\def\L{{\bf L}}
\def\E{{\mathbb E}}
\def\F{{\bf F}}
\def\P{{\mathbb P}}
\def\N{{\mathbb N}}
\def\eps{\varepsilon}
\def\wh{\widehat}
\def\pf{\noindent{\bf Proof.} }

\title{\Large \bf Boundary Harnack Principle for Subordinate Brownian Motions}
\author{Panki Kim\thanks{The research of this author is supported
by the Korea Research Foundation Grant funded by the Korean Government
(MOEHRD, Basic Research Promotion Fund) (KRF-2007-331-C00037).}\\
Department of Mathematics\\
Seoul National University\\
Seoul 151-742, Republic of Korea\\
Email: pkim@snu.ac.kr, \smallskip \\
\smallskip \\
Renming Song\\
Department of Mathematics\\
University of Illinois \\
Urbana, IL 61801, USA\\
Email: rsong@math.uiuc.edu \smallskip \\
and \smallskip \\
Zoran Vondra\v{c}ek\thanks{The research of this author is supported in part by
MZOS grant 037-0372790-2801 of the Republic of Croatia.}\\
Department of Mathematics\\
University of Zagreb \\
Zagreb, Croatia\\
Email: vondra@math.hr }

\date{ }
\maketitle

\begin{abstract}
We establish a boundary Harnack principle for a large class of subordinate
Brownian motion, including mixtures of symmetric stable processes,
in bounded $\kappa$-fat open set (disconnected analogue of John domains).
As an application of the boundary Harnack principle,  we identify the
Martin boundary and the minimal Martin boundary of bounded $\kappa$-fat
open sets with respect to these processes
with their Euclidean boundary.
\end{abstract}

\noindent {\bf AMS 2000 Mathematics Subject Classification}: Primary 60J45,
Secondary 60J25, 60J51.

\noindent {\bf Keywords and phrases:}
Green functions, Poisson kernels, subordinator, subordinate
Brownian motion, Bernstein functions, complete Bernstein functions,
 symmetric stable processes, mixture of symmetric stable processes,
harmonic functions, Harnack inequality, boundary Harnack principle,
Martin boundary.

\section{Introduction}

The boundary Harnack principle
for nonnegative classical harmonic functions is a very deep result in
potential theory and has very important applications in probability
and potential theory.

In \cite{B} Bogdan showed that the boundary Harnack principle is valid
in bounded Lipschitz domains for nonnegative harmonic functions of
rotationally invariant stable processes and
then in \cite{SW} Song and Wu extended
the boundary Harnack principle for rotationally invariant
stable processes to bounded $\kappa$-fat open set.
Subsequently Bogdan-Stos-Sztonyk \cite{BSS} and Sztonyk \cite{Sz2}
extended the boundary Harnack principle to symmetric (not necessarily
rotationally invariant) stable processes. In a recent paper \cite{BKK},
Bogdan, Kulczycki and Kwasnicki proved a version of the boundary Harnack
inequality for nonnegative harmonic functions of
rotationally invariant stable processes in arbitrary open sets.

By using some perturbation methods, the boundary Harnack principle
has been generalized to some classes of  rotationally invariant
L\`evy processes including relativistic stable processes and
truncated stable processes. These processes can be regarded as
perturbations of rotationally invariant
stable processes and their Green functions on bounded smooth domains
are comparable to their counterparts for rotationally invariant
 stable processes
(see \cite{CS3}, \cite{R}, \cite{GR}, \cite{KL}, \cite{KS1} and \cite{KS2}).
This comparison of Green functions played a crucial role in the arguments
of \cite{GR} ,  \cite{KS1} and \cite{KS2}.

In this paper, we will show that, under minimal conditions,
the boundary Harnack principle is
valid for subordinate Brownian motions with characteristic exponents
of the form $\Phi(\xi)=|\xi|^{\alpha}\ell(|\xi|^2)$ for some $\alpha\in
(0, 2)$ and some positive function $\ell$ which is slowly varying at
$\infty$. Examples of this class of subordinate Brownian motions include,
among others, relativistic stable processes and mixtures of rotationally
invariant stable processes.
The  Green functions of subordinate Brownian motions considered here
behave like $\frac{c}{|x|^{d-\alpha}\ell(|x|^{-2})}$ near the origin.
So these subordinate Brownian motions
can not be regarded as perturbations of rotationally invariant
stable processes in general
and their Green functions in bounded smooth domains are
not comparable to their counterparts for rotationally invariant
stable processes.

Our proof of the boundary Harnack principle will be similar to the
arguments in \cite{B} and \cite{SW} for rotationally invariant
stable processes.  One of the key ingredients is a sharp upper
bound for the expected exit time from a ball which, in the case of stable
processes, follows easily from the explicit formula for the Green
function of a ball. However, in the present case, the desired upper
bound is pretty difficult to establish. We rely on
the fluctuation theory for real-valued L\'evy processes to accomplish
this.

The organization of this paper is as follows. In Section 2 we use
the fluctuation theory for real-valued L\'evy processes to establish
a nice upper bound on the expected exit time from an interval for
a one-dimensional subordinate Brownian motion. In Section 3, we use
the results of Section 2 to establish the desired upper bound
on the expected exit time from a ball for a multidimensional
subordinate Brownian motion and an upper bound on the Poisson
kernel of a ball. The proof of the boundary Harnack principle
is given in Section 4 and in the last section we apply our boundary
Harnack principle to study the Martin boundary with respect to
subordinate Brownian motions.

In this paper we will use the following convention: the values of
the constants $r_1, r_2, \dots$ will remain the same throughout
this paper, while the values of the constants $c_1, c_2, \dots$
or $C, C_1, C_2, \dots$
might change from one appearance to another.
The dependence of the constants on the dimension, the index $\alpha$
and the slowly varying function will not be mentioned explicitly,
while the dependence of the constants on other quantities will
be expressed using $c(\cdot)$ with the arguments representing
the quantities the constant depends on.
In this paper, we use ``$:=$" to denote a definition, which
is  read as ``is defined to be".
$f(t) \sim g(t)$, $t \to 0$ ($f(t) \sim g(t)$, $t \to \infty$,
respectively) means $
\lim_{t \to 0} f(t)/g(t) = 1$ ($\lim_{t \to \infty} f(t)/g(t) = 1$,
respectively).

\section{Some Results on One-dimensional Subordinate Brownian Motion}

Suppose that $W=(W_t: t\ge 0)$ is a one-dimensional Brownian motion with
$$
\E\left[e^{i\xi(W_t-W_0)}\right]=e^{-t\xi^2}, \qquad
\forall \xi\in \R,~ t>0\, ,
$$
and $S=(S_t: t\ge 0)$ is a subordinator
(a non-negative increasing L\'evy process)
independent of $W$ and with
Laplace exponent $\phi$, that is
$$
\E\left[e^{-\lambda S_t}\right]=e^{-t\phi(\lambda)}, \qquad
\forall t, \lambda>0.
$$
A $C^{\infty}$ function $g : (0, \infty) \to [0, \infty)$ is called a
Bernstein function if $(-1)^n D^ng \le 0$ for every positive integer $n$.
A Bernstein function $g$ can be written in the following form
$$
g(\lambda)=a+b\lambda+\int^{\infty}_0(1-e^{-\lambda t})\mu(dt)
$$
where $a, b\ge 0$ and $\mu$ is a measure on $(0, \infty)$ with
$\int^{\infty}_0(1\wedge t)\mu(dt)<\infty$. $\mu$ is called the L\`evy
measure of $g$. It is well known that a function $g$ is the
Laplace exponent of a (non-killed) subordinator if and only if $g$ is a Bernstein function
with $\lim_{\lambda\to 0}g(\lambda)=0$.
A Bernstein function $g$ is called a complete Bernstein function if its
L\'evy measure $\mu$ has a completely monotone density with respect
to the Lebesgue measure. For details on examples and properties of
complete Bernstein functions, one can see \cite{Jac}, \cite{Sch} or
\cite{SV2}. One important property we are going to use in this paper
is that $f$ is complete Bernstein is equivalent to that $\lambda/f(\lambda)$
is complete Bernstein.
Throughout this paper we will assume that $\phi$ is a complete
Bernstein function such that
\begin{equation}\label{e:asumptiononle}
\phi(\lambda)=\lambda^{\alpha/2}\ell(\lambda)
\end{equation}
for some $\alpha\in (0, 2)$ and some positive function $\ell$ which
is slowly varying at $\infty$. For concepts and results related to the slowly varying function, we refer our
readers to \cite{BGT}.

Using Corollary 2.3 of \cite{SV3} or Theorem 2.3 of \cite{RSV}
we know that the potential measure of $S$ has a decreasing density $u$.
By using the Tauberian theorem (Theorem 1.7.1 in \cite{BGT}) and
the monotone density
theorem (Theorem 1.7.2 in \cite{BGT}), one can easily check that
\begin{equation}\label{e:abofpdofSat0}
u(t)\sim\frac{t^{\alpha/2-1}}{\Gamma(\alpha/2)}\frac1{\ell(t^{-1})},
\qquad t\to 0.
\end{equation}
Let $\mu(t)$ be the density of the  L\'evy measure of $\phi$.
It follows from
Proposition 2.23 of \cite{SV2} that
\begin{equation}\label{e:abofldofSat0}
\mu(t)\sim\frac{\alpha}{2\Gamma(1-\alpha/2)}\frac{\ell(t^{-1})}{
t^{1+\alpha/2}},
\qquad t\to 0.
\end{equation}

The subordinate Brownian motion $X=(X_t:t\ge 0)$ defined by
$X_t=W_{S_t}$ is a symmetric L\'evy process with the characteristic
exponent
$$
\Phi(\theta)=\phi(\theta^2)=|\theta|^{\alpha}\ell(\theta^2),
\qquad \forall \theta\in \R.
$$
Let $\chi$ be the Laplace exponent of the ladder height process
of $X$.
(For the definition of the ladder height process and its basic properties,
we refer our readers to Chapter 6 of \cite{Ber}.)
Then it follows from Corollary 9.7 of \cite{Fris} that
\begin{equation}\label{e:formula4leoflh}
\chi(\lambda)=
\exp\left(\frac1\pi\int^{\infty}_0\frac{\log(\Phi(\lambda\theta))}
{1+\theta^2}d\theta
\right)
=\exp\left(\frac1\pi\int^{\infty}_0\frac{\log(
|\theta|^{\alpha}\lambda^{\alpha}\ell(\theta^2\lambda^2))}{1+\theta^2}d\theta
\right), \quad \forall \lambda>0.
\end{equation}
Under our assumptions, we have the following result.

\begin{prop}\label{p:lhisspecial}
The Laplace exponent $\chi$ of the ladder height process of $X$ is
a special Bernstein function.
i.e., $\lambda/ \chi(\lambda)$ is also a Bernstein function.
\end{prop}

\pf Define $\psi(\lambda)=\lambda/\phi(\lambda)$. Let $T$ be a subordinator
independent of $W$ and with Laplace exponent $\psi$ and let $Y=(Y_t:t\ge 0)$
be the subordinate Brownian motion defined by $Y_t=W_{T_t}$. Let $\Psi$ be
the characteristic exponent of $Y$. Then
$$
\Phi(\theta)\Psi(\theta)=\phi(\theta^2)\psi(\theta^2)=\theta^2,
\qquad \forall \theta\in \R.
$$
Let $\rho$ be the Laplace exponent of the ladder height process of $Y$. Then
by (\ref{e:formula4leoflh}) we have
\begin{eqnarray*}
\chi(\lambda)\rho(\lambda)&=& \exp\left(\frac{1}{\pi}
\int_0^{\infty} \frac{\log(\Phi(\theta\lambda))
+\log(\Psi(\theta\lambda))}{1+\theta^2}\, d\theta\right)\\
&=&\exp\left(\frac{1}{\pi} \int_0^{\infty}
\frac{\log(\Phi(\theta\lambda)\Psi(\theta\lambda))}{1+\theta^2}\,
d\theta\right)\\
&=&\exp\left(\frac{1}{\pi} \int_0^{\infty}
\frac{\log(\theta^2\lambda^2)}{1+\theta^2}\, d\theta\right)=\lambda.
\end{eqnarray*}
Thus $\chi$ is a special Bernstein function.
\qed

\begin{prop}\label{p:abofkappaatinfty}
If there are $M>1$, $\delta\in (0, 1)$ and
a nonnegative integrable function $f$ on $(0, \delta)$ such that
\begin{equation}\label{e:cdtogetabofkappa}
\left|\log\left(\frac{\ell(\lambda^2\theta^2)}{\ell(\lambda^2)}\right)\right|
\le f(\theta), \qquad \forall (\theta, \lambda) \in (0, \delta)
\times (M, \infty),
\end{equation}
then
\begin{equation}\label{e:abofkappaatinfty}
\lim_{\lambda\to\infty}\frac{\chi(\lambda)}{\lambda^{\alpha/2}(
\ell(\lambda^2))^{1/2}}=1.
\end{equation}
\end{prop}

\pf
Using the identity
$$
\lambda^{\beta/2}=\exp\left(\frac{1}{\pi}\int_0^{\infty}
\frac{\log(\theta^{\beta}\lambda^{\beta})}
{1+\theta^2}\, d\theta\right),\, \qquad \forall \lambda, \beta>0,
$$
we get easily from (\ref{e:formula4leoflh}) that
\begin{eqnarray*}
\chi(\lambda)&=&\lambda^{\alpha/2}
\exp\left(\frac1\pi\int^{\infty}_0\frac{\log(\ell(\lambda^2\theta^2))}
{1+\theta^2}d\theta\right)\\
&=&\lambda^{\alpha/2}(\ell(\lambda^2))^{1/2}
\exp\left(\frac1\pi\int^{\infty}_0\log\left(\frac{\ell(\lambda^2\theta^2)}
{\ell(\lambda^2)}\right)\frac1
{1+\theta^2}d\theta\right).
\end{eqnarray*}
By Potter's Theorem (Theorem 1.5.6 (1) in \cite{BGT}),
there exists $\lambda_0 >1$ such that
$$
\left|\log\left(\frac{\ell(\lambda^2\theta^2)}
{\ell(\lambda^2)}\right)\right|\frac1
{1+\theta^2} \, \le\, 2 \,\frac {\log \theta}
{1+\theta^2}, \qquad  \forall (\theta, \lambda) \in
[1, \infty) \times [\lambda_0, \infty).
$$
Thus by using the dominated convergence theorem in the first integral below,
the uniform convergence theorem
(Theorem 1.2.1 in \cite{BGT}) in the second integral, and the assumption
(\ref{e:cdtogetabofkappa}) in the third integral, we have
$$
\lim_{\lambda \to \infty} \int^{\infty}_{0}\log\left
(\frac{\ell(\lambda^2\theta^2)}
{\ell(\lambda^2)}\right)\frac1
{1+\theta^2}\, d\theta=\lim_{\lambda \to \infty}
\left(\int^{\infty}_{1} + \int^{1}_{\delta}+\int_0^{\delta}\right)\log
\left(\frac{\ell(\lambda^2\theta^2)}
{\ell(\lambda^2)}\right)\frac1
{1+\theta^2}\, d\theta=0.
$$
\qed

In the case $\phi(\lambda)=\lambda^{\alpha/2}$ for some $\alpha\in (0, 2)$,
the assumption of
the proposition above is trivially satisfied. Now we give some
other examples.

\begin{example}\label{relstablesub}
{\rm Suppose that $\alpha\in (0, 2)$ and define
$$
\phi(\lambda)=(\lambda+1)^{\alpha/2}-1.
$$
Then $\phi$ is a complete Bernstein function which can be written
as $\phi(\lambda)=\lambda^{\alpha/2}\ell(\lambda)$ with
$$
\ell(\lambda)=\frac{(\lambda+1)^{\alpha/2}-1}{\lambda^{\alpha/2}}.
$$
Using elementary analysis one can easily check that there is a
nonnegative integrable function $f$ on $(0, 1)$ such
that (\ref{e:cdtogetabofkappa}) is satisfied.}
\end{example}

\begin{example}\label{stablemix}
{\rm Suppose $0<\beta<\alpha<2$ and define
$$
\phi(\lambda)=\lambda^{\alpha/2}+\lambda^{\beta/2}.
$$
Then $\phi$ is a complete Bernstein function which can be written
as $\phi(\lambda)=\lambda^{\alpha/2}\ell(\lambda)$ with
$$
\ell(\lambda)=1+\lambda^{(\beta-\alpha)/2}.
$$
Using elementary analysis one can easily check that there is a
nonnegative integrable function $f$ on $(0, 1)$ such
that (\ref{e:cdtogetabofkappa}) is satisfied.}
\end{example}

\begin{example}\label{stablewlogc1}
{\rm Suppose that $\alpha\in (0, 2)$ and $\beta\in (0, 2-\alpha)$.
Define
$$
\phi(\lambda)=\lambda^{\alpha/2}(\log(1+\lambda))^{\beta/2}.
$$
By using the facts that $\lambda$ and $\log(1+\lambda)$ are complete
Bernstein functions and properties of complete Bernstein functions
(see \cite{SV2}), one can easily check that $\phi$ is a complete Bernstein
function. $\phi$ can be written
as $\phi(\lambda)=\lambda^{\alpha/2}\ell(\lambda)$ with
$$
\ell(\lambda)=(\log(1+\lambda))^{\beta/2}.
$$
To check that there is a
nonnegative integrable function $f$ on $(0, 1)$ such
that (\ref{e:cdtogetabofkappa}) is satisfied, we only
need to bound the function
$$
\left|\log\left(\frac{\log(1+\lambda^2\theta^2)}{\log(1+\lambda^2)}
\right)\right|
$$
for large $\lambda$ and small $\theta$. We will consider two cases
separately. Fix an $M>1$ and a $\theta<1$.
\begin{description}
\item{(1)} $\lambda\ge M$, $\theta<1$ and $\lambda>1/\theta$. In
this case, by using the fact that for any $a>0$ the function
$x\mapsto \frac{x}{x-a}$ is decreasing on $(a, \infty)$, we get that
\begin{eqnarray*}
\left|\log\left(\frac{\log(1+\lambda^2\theta^2)}{\log(1+\lambda^2)}
\right)\right|&=&
\log\left(\frac{\log(1+\lambda^2)}{\log(1+\lambda^2\theta^2)}\right)\\
&\le&\log\left(\frac{\log(1+\lambda^2)}{\log(\theta^2)+
\log(1+\lambda^2)}\right)\\
&\le&\log\left(\frac{\log(1+\theta^{-2})}{\log(\theta^2)+
\log(1+\theta^{-2})}\right)\\
&=&\log\left(\frac{\log(1+\theta^{2})-\log(\theta^2)}{
\log(1+\theta^{2})}\right).
\end{eqnarray*}

\item{(2)} $\lambda\ge M$, $\theta<1$ and $\lambda\le1/\theta$. In
this case we have
\begin{eqnarray*}
\left|\log\left(\frac{\log(1+\lambda^2\theta^2)}{\log(1+\lambda^2)}
\right)\right|&=&
\log\left(\frac{\log(1+\lambda^2)}{\log(1+\lambda^2\theta^2)}\right)\\
&\le& \left(\frac{\log(1+\lambda^2)}{\log(1+M^2\theta^2)}\right)\\
&\le& \left(\frac{\log(1+\theta^{-2})}{\log(1+M^2\theta^2)}\right).
\end{eqnarray*}
\end{description}
Combining the results above one can easily check that
there is a
nonnegative integrable function $f$ on $(0, 1)$ such
that (\ref{e:cdtogetabofkappa}) is satisfied.
}
\end{example}

\begin{example}\label{stablewlogc2}
{\rm Suppose that $\alpha\in (0, 2)$ and $\beta\in (0, \alpha)$.
Define
$$
\phi(\lambda)=\lambda^{\alpha/2}(\log(1+\lambda))^{-\beta/2}.
$$
By using the facts that $\lambda$ and $\log(1+\lambda)$ are complete
Bernstein functions and properties of complete Bernstein functions
(see \cite{SV2}), one can easily check that $\phi$ is a complete Bernstein
function. $\phi$ can be written
as $\phi(\lambda)=\lambda^{\alpha/2}\ell(\lambda)$ with
$$
\ell(\lambda)=(\log(1+\lambda))^{-\beta/2}.
$$
Similarly to the example above, once can use elementary analysis
to check that there is a
nonnegative integrable function $f$ on $(0, 1)$ such
that (\ref{e:cdtogetabofkappa}) is satisfied.
}
\end{example}

The method used to construct the complete Bernstein functions can
be used to construct a whole class of complete Bernstein functions
satisfying the assumptions of this paper. For instance, one can
check that, for $\alpha\in (0, 2)$, $\beta\in (0, 2-\alpha)$, functions
like $\lambda^{\alpha/2}(\log(1+\log(1+\lambda)))^{\beta/2},
\lambda^{\alpha/2}(\log(1+\log(1+\log(1+\lambda))))^{\beta/2},
\dots$ are complete Bernstein functions satisfying the assumptions
of this paper. Similarly, for any $\alpha\in (0, 2)$,
$\beta\in (0, \alpha)$, functions
like $\lambda^{\alpha/2}(\log(1+\log(1+\lambda)))^{-\beta/2},
\lambda^{\alpha/2}(\log(1+\log(1+\log(1+\lambda))))^{-\beta/2},
\dots$ are complete Bernstein functions satisfying the assumptions
of this paper.

In the remainder of this section we will always assume that the assumption
of Proposition \ref{p:abofkappaatinfty} is satisfied.
It follows from Propositions \ref{p:lhisspecial} and \ref{p:abofkappaatinfty}
above and Corollary 2.3 of \cite{SV3}
that the potential measure $V$ of the ladder height process of $X$ has a
decreasing density $v$. Since $X$ is symmetric, we know
that the potential measure $\hat V$ of the dual ladder height process
is equal to $V$.

In light of Proposition \ref{p:abofkappaatinfty}, one can easily
apply the Tauberian theorem (Theorem 1.7.1 in \cite{BGT}) and
the monotone density
theorem (Theorem 1.7.2 in \cite{BGT})
to get the following result.

\begin{prop}\label{p: abofgf4lhpat0}
As $x\to 0$, we have
\begin{eqnarray*}
V((0, x))&\sim&\frac{x^{\alpha/2}}{\Gamma(1+\alpha/2)
(\ell(x^{-2}))^{1/2}},\\
v(x)&\sim&\frac{x^{\alpha/2-1}}{\Gamma(\alpha/2)
(\ell(x^{-2}))^{1/2}}.
\end{eqnarray*}
\end{prop}

\pf We omit the details.
\qed

It follows from Proposition \ref{p:abofkappaatinfty} above and
Lemma 7.10 of \cite{Ky} that the process $X$ does not creep upwards.
Since $X$ is symmetric, we know that $X$ also does not creep
downwards. Thus if, for any $a\in \R$, we define
$$
\tau_a=\inf\{t>0: X_t<a\}, \quad \sigma_a=\inf\{t>0: X_t\le a\},
$$
then we have
\begin{equation}\label{e:firstexittime}
\P_x(\tau_a=\sigma_a)=1, \quad x>a.
\end{equation}

Let $G^{(0, \infty)}(x, y)$ be the Green function of $X^{(0, \infty)}$,
the process obtained by killing $X$ upon exiting
from $(0, \infty)$. Then we have the following result.

\begin{prop}\label{p:Greenf4kpXonhalfline} For any $x, y>0$ we have
$$
G^{(0, \infty)}(x, y)=\left\{\begin{array}{ll}
\int^x_0v(z)v(y+z-x)dz, & x\le y,\\
\int^x_{x-y}v(z)v(y+z-x)dz, & x>y.
\end{array}\right.
$$
\end{prop}

\pf By using (\ref{e:firstexittime}) above and Theorem 20 on page
176 of \cite{Ber} we get that for any nonnegative function on $f$
on $(0, \infty)$,
\begin{equation}\label{e:e1inpfoformforgf}
\E_x\left[ \int_0^{\infty} f(X^{(0, \infty)}_t)\, dt\right]=k
\int^{\infty}_0
\int^x_0v(z)f(x+z-y)v(y) dz dy\, ,
\end{equation}
where $k$ is the constant depending on the normalization of the local time
of the process $X$ reflected at its supremum.
We choose $k=1$. Then
\begin{eqnarray}\label{e:e2inpfoformforgf}
&&\E_x\left[ \int_0^{\infty} f(X^{(0, \infty)}_t)\, dt\right]
\,=\,\int_0^{\infty} \, v(y)\int_0^x \, v(z)f(x+y-z) dz dy\nonumber\\
&&=\int_0^x \, v(z)\int_0^{\infty}  v(y)f(x+y-z) dy dz
\,=\,\int_0^x \, v(z)\int_{x-z}^{\infty}\, v(w+z-x)f(w) dw dz\nonumber\\
&&=\int_0^x f(w) \int_{x-w}^x \,  v(z)v(w+z-x) dzdw+
\int_x^{\infty} f(w) \int_0^x \, v(z)v(w+z-x)  dzdw\, .
\end{eqnarray}
On the other hand,
\begin{eqnarray}\label{e:e3inpfoformforgf}
&&\E_x\left[\int_0^{\infty} f(X^{(0, \infty)}_t)\, dt\right]
=\int_0^{\infty}G^{(0, \infty)}(x,w)f(w)\, dw\nonumber \\
&&=\int_0^x G^{(0, \infty)}(x,w)f(w)\, dw +
\int_x^{\infty} G^{(0, \infty)}(x,w)f(w)\, dw\, .
\end{eqnarray}
By comparing (\ref{e:e2inpfoformforgf}) and (\ref{e:e3inpfoformforgf})
we arrive at our desired conclusion.
\qed

For any $r>0$, let $G^{(0, r)}$ be the Green function of $X^{(0, r)}$,
the process obtained by killing $X$ upon exiting from $(0, r)$.
Then we have the following result.

\begin{prop}\label{p:upbdongfofkpinfiniteinterval}
For any $R>0$, there exists $C=C(R)>0$ such that
$$
\int^r_0G^{(0, r)}(x, y)dy\le C\frac{r^{\alpha/2}}{(\ell(r^{-2}))^{1/2}}
\frac{x^{\alpha/2}}{(\ell(x^{-2}))^{1/2}}, \qquad x\in (0, r),~ r\in (0, R).
$$
\end{prop}

\pf For any $x\in (0, r)$, we have
\begin{eqnarray*}
&&\int^r_0G^{(0, r)}(x, y)dy
\le \int^r_0G^{(0, \infty)}(x, y)dy\\
&&=\int^x_0\int^x_{x-y}v(z)v(y+z-x)dzdy+
\int^r_x\int^x_0v(z)v(y+z-x)dzdy\\
&&=\int^x_0v(z)\int^x_{x-z}v(y+z-x)dydz
+\int^x_0v(z)\int^r_xv(y+z-x)dydz\,
\le\, 2\,V((0, r))\,V((0, x)).
\end{eqnarray*}
Now the desired conclusion follows easily from Proposition
\ref{p: abofgf4lhpat0} and the continuity of $V((0, x))$
and $x^{\alpha/2}/(\ell(x^{-2}))^{1/2}$.

\qed

As a consequence of the result above, we immediately get the following.

\begin{prop}\label{p:upbdongfofkpinfiniteinterval2}
For any $R>0$, there exists $C=C(R)>0$ such that
$$
\int^r_0G^{(0, r)}(x, y)dy\le C\frac{r^{\alpha/2}}{(\ell(r^{-2}))^{1/2}}
\left(\frac{x^{\alpha/2}}{(\ell(x^{-2}))^{1/2}}\wedge
\frac{(r-x)^{\alpha/2}}{(\ell((r-x)^{-2}))^{1/2}}\right),
\quad x\in (0, r), ~r\in (0, R).
$$
\end{prop}

\section{Key estimates on Multi-dimensional Subordinate Brownian Motions}

In the remainder of this paper we will always assume that
$d\ge 2$ and that $\alpha \in (0, 2)$. From now on we will assume that
$B=(B_t:t\ge 0)$ is a Brownian motion on $\R^d$ with
$$
\E\left[e^{i\xi\cdot(B_t-B_0)}\right]=
e^{-t|\xi|^2}, \qquad \forall \xi \in \R^d, t>0.
$$
Suppose that $S=(S_t: t\ge 0)$ is a subordinator independent of $B$
and that its Laplace exponent $\phi$ is a complete Bernstein function
satisfying all the assumption of the previous section. More precisely
we assume that there is a positive function $\ell$ on $(0, \infty)$
which is slowly varying at $\infty$ such that $\phi(\lambda)=
\lambda^{\alpha/2}\ell(\lambda)$ for all $\lambda>0$
and that there is a nonnegative integrable function $f$
on $(0, \delta)$ for some $\delta>0$
such that (\ref{e:cdtogetabofkappa}) holds.
As in the previous section, we will use $u(t)$ and $\mu(t)$
to denote the potential density and L\'evy density of $S$ respectively.

In the sequel, we will use $X=(X_t: t\ge 0)$ to denote
the subordinate Brownian motion defined by $X_t=B_{S_t}$. Then it is easy
to check that when $d\ge 3$ the process $X$ is transient. In the case
of $d=2$, we will always assume the following:

{\bf A1}. {\it The potential density $u$ of
$S$ satisfies the following assumption:
\begin{equation}\label{e:abatinftyofptofS}
u(t)\sim ct^{\gamma-1}, \qquad t\to\infty
\end{equation}
for some constants $c>0$
and $\gamma < 1$.}

Under this assumption the process $X$ is
also transient for $d=2$.

We will use $G(x, y)=G(x-y)$
to denote the Green function of $X$.
 The Green function $G$ of $X$ is given by the following formula
$$
G(x)=\int^{\infty}_0(4\pi t)^{-d/2}e^{-|x|^2/(4t)}u(t)dt, \qquad x\in \R^d.
$$
Using this formula, we can easily see that $G$ is
radially decreasing and continuous
in $\R^d\setminus \{0\}$.

In order to get the asymptotic behavior of $G$ near the origin, we
need some additional assumption on the slowly varying function $\ell$.
For any $y, t, \xi>0$, define
$$
\Lambda_{\ell,\xi}(y,t):=\left\{\begin{array}{cl}
\frac{\ell(1/y)}{\ell(4t/y)}, &y< \frac{t}{\xi}\, ,\\
0, &y\ge \frac{t}{\xi}\, .
\end{array}\right.
$$
We will always assume that

{\bf A2}. {\it There is a $\xi >0$ such that
$$
\Lambda_{\ell, \xi}(y, t)\le g(t), \qquad \forall y, t>0,
$$
for some positive function $g$ on $(0, \infty)$ with
$$
\int^{\infty}_0t^{(d-\alpha)/2-1}e^{-t}g(t)dt<\infty.
$$}

It is easy to check (see the proofs of Theorem 3.6 and
Theorem 3.11 in \cite{SV2}) that for the subordinators corresponding
to Examples \ref{relstablesub}--\ref{stablewlogc2}, {\bf A1}  and
{\bf A2} are satisfied.

Under these assumptions we have the following.

\begin{thm}\label{t:abofgfatorigin}
The Green function $G$ of $X$ satisfies the
following
$$
G(x)\sim \frac{\alpha\Gamma((d-\alpha)/2)}{2^{\alpha+1}\pi^{d/2}\Gamma(1+
\alpha/2)}\frac1{|x|^{d-\alpha}\ell(|x|^{-2})},
\qquad |x|\to 0.
$$
\end{thm}

\pf This follows easily from {\bf A1}-{\bf A2}, (\ref{e:abofpdofSat0})
above and Lemma
3.3 of \cite{SV2}. We omit the details.
\qed

Let $J$ be the jumping function of $X$,
then
$$
J(x)=\int^{\infty}_0(4\pi t)^{-d/2}e^{-|x|^2/(4t)}\mu(t)dt, \qquad x\in \R^d.
$$
Thus $J(x)=j(|x|)$ with
$$
j(r)=\int^{\infty}_0(4\pi t)^{-d/2}e^{-r^2/(4t)}\mu(t)dt, \qquad r>0.
$$
It is easy to see that $j$ is continuous in $(0, \infty)$.
Since $t\mapsto \mu(t)$ is decreasing, the function $r\mapsto
j(r)$ is decreasing on $(0, \infty)$.
In order to get the asymptotic behavior of $j$ near the origin, we
need some additional assumption on the slowly varying function $\ell$.
For any $y, t, \xi>0$, define
$$
\Upsilon_{\ell,\xi}(y,t):=\left\{\begin{array}{cl}
\frac{\ell(4t/y)}{\ell(1/y)}, &y< \frac{t}{\xi}\, ,\\
0, &y\ge \frac{t}{\xi}\, .
\end{array}\right.
$$
We will always assume that

{\bf A3}. {\it There is a $\xi >0$ such that
$$
\Upsilon_{\ell, \xi}(y, t)\le h(t), \qquad \forall y, t>0
$$
for some positive function $h$ on $(0, \infty)$ with
$$
\int^{\infty}_0t^{(d+\alpha)/2-1}e^{-t}h(t)dt<\infty.
$$}

It is easy to check (see the proofs of Theorem 3.6 and
Theorem 3.11 in \cite{SV2}) that for the subordinators corresponding
to Examples \ref{relstablesub}--\ref{stablewlogc2}, {\bf A3} is satisfied.

\begin{thm}\label{t:abofjfatorigin}
The function $j$ satisfies the following
$$
j(r)\sim \frac{\alpha\Gamma((d+\alpha)/2)}{2^{1-\alpha}\pi^{d/2}\Gamma(1-
\alpha/2)}\frac{\ell(r^{-2})}{r^{d+\alpha}},
\qquad r\to 0.
$$
\end{thm}

\pf This follows easily from {\bf A1}, {\bf A3}, (\ref{e:abofldofSat0})
above and Lemma 3.3 of \cite{SV2}. We omit the details.
\qed

For any open set $D$, we use $\tau_D$ to denote the first exit
time from $D$, i.e., $\tau_D=\inf\{t>0: \, X_t\notin D\}$.
Given  an open set $D\subset \R^d$, we define
$X^D_t(\omega)=X_t(\omega)$ if $t< \tau_D(\omega)$ and
$X^D_t(\omega)=\partial$ if $t\geq  \tau_D(\omega)$, where
$\partial$ is a cemetery state. We now recall the definition
of harmonic functions with respect to $X$.

\begin{defn}\label{def:har1}
Let $D$ be an open subset of $\R^d$.
A function $u$ defined on $\R^d$ is said to be

\begin{description}
\item{(1)}  harmonic in $D$ with respect to $X$ if
$$
\E_x\left[|u(X_{\tau_{B}})|\right] <\infty
\quad \hbox{ and } \quad
u(x)= \E_x\left[u(X_{\tau_{B}})\right],
\qquad x\in B,
$$
for every open set $B$ whose closure is a compact
subset of $D$;

\item{(2})
regular harmonic in $D$ with respect to $X$ if it is harmonic in $D$
with respect to $X$ and
for each $x \in D$,
$$
u(x)= \E_x\left[u(X_{\tau_{D}})\right];
$$

\item{(3}) harmonic for $X^D$ if it is harmonic for $X$ in $D$ and
vanishes outside $D$.

\end{description}
\end{defn}

In order for a scale invariant Harnack inequality to hold, we need
to assume some additional conditions on the L\'evy density $\mu$
of $S$. We will always assume that

{\bf A4}. {\it The L\'evy density $\mu$ of $S$ satisfies the following
conditions: there exists $C_1>0$ such that
$$
\mu(t)\le C_1\mu(t+1), \qquad \forall t>1.
$$
}

It follows from \eqref{e:abofldofSat0} that for any $M>0$ there
exists $C_2>0$ such that
$$
\mu(t)\le C_2\mu(2t), \qquad \forall t\in (0, M).
$$

Using {\bf A4}, the display above and repeating the proof of
Lemma 4.2 of \cite{RSV} we get that

\begin{description}
\item{(1)} For any $M>0$, there exists $C_3>0$ such that
\begin{equation}\label{H:1}
j(r)\le C_3 j(2r), \qquad \forall r\in (0, M).
\end{equation}
\item{(2)} There exists $C_4>0$ such that
$$
j(t)\le C_4 j(r+1), \qquad \forall r>1.
$$
\end{description}

It is easy to check (see \cite{SV2}) that for the subordinators corresponding
to Examples \ref{relstablesub}--\ref{stablewlogc2}, {\bf A4} is satisfied.
Therefore by Theorem 4.14 of \cite{SV2} (see also \cite{RSV}) we have
the following Harnack inequality.

\begin{thm}[Harnack inequality]\label{T:Har}
There exist $r_1 \in  (0,1)$ and $C>0$  such that for every $r\in (0,r_1)$,
every $x_0\in \R^d$, and every nonnegative function $u$ on $\R^d$
which is harmonic
in $B(x_0,r)$ with respect to $X$, we have
$$
\sup_{y \in B(x_0, r/2)}u(y) \le C \inf _{y \in B(x_0, r/2)}u(y).
$$
\end{thm}

For any bounded open set $D$ in $\R^d$, we will use $G_D(x,y)$
to denote the Green function of $X^D$. Using the continuity
and the radial decreasing property of $G$, we can easily check that
$G_D$ is continuous in $(D\times D)\setminus\{(x, x): x\in D\}$.

\begin{prop}\label{p:Green}
For any $R>0$, there exists $C=C(R)>0$ such that for every open
subset $D$ with diam$(D)\le R$,
\begin{equation}\label{G2}
G_D(x,y) \,\le\, G(x,y) \,\le\, C \,\frac{1}
{\ell(|x-y|^{-2})|x-y|^{d-\alpha}}, \qquad\forall
(x,y) \in D \times D.
\end{equation}
\end{prop}

\pf
The results of this proposition are immediate consequences of Theorem
\ref{t:abofgfatorigin} and the continuity and positivity of
$\ell(r^{-2}) r^{d-\alpha}$ on $(0, \infty)$.
\qed

The idea of the proof of the next lemma comes from \cite{Sz2}.
\begin{lemma}\label{l:tau}
For any $R>0$, there exists $C=C(R)>0$ such that for every $r\in (0, R)$
and $x_0 \in \R^d$,
$$\E_x[\tau_{B(x_0,r)}]\,\le\, C\, \frac{r^{\alpha/2}}
{(\ell(r^{-2}))^{1/2}}\frac{(r-|x-x_0|)^{\alpha/2}}
{(\ell((r-|x-x_0|)^{-2}))^{1/2}}, \qquad x\in B(x_0, r).
$$
\end{lemma}

\pf Without loss of generality, we may assume that $x_0=0$.
For $x\neq 0$, put $Z_t=\frac{X_t\cdot x}{|x|}$. Then $Z_t$ is a
L\'evy process on $\R$ with
$$
\E(e^{i\theta Z_t})=\E(e^{i\theta\frac{x}{|x|}\cdot X_t})
=e^{-t|\theta|^{\alpha}\ell(\theta^2)}, \qquad \theta\in \R.
$$
Thus $Z_t$ is the type of one-dimensional subordinate Brownian motion
we studied
in the previous section. It is easy to see that, if $X_t\in B(0, r)$,
then $|Z_t|<r$, hence
$$
\E_x[\tau_{B(0, r)}]\le \E_{|x|}[\tilde \tau],
$$
where $\tilde \tau=\inf\{t>0: |Z_t|\ge r\}$. Now the desired conclusion
follows easily from Proposition \ref{p:upbdongfofkpinfiniteinterval2}.
\qed

\begin{lemma}\label{l:tau2}
There exist $r_2 \in (0,r_1]$ and $C>0$ such that for every positive
$r \le r_2$ and $x_0 \in \R^d$,
$$\E_{x_0}[\tau_{B(x_0,r)}]\,\ge\, C\, \frac{r^\alpha}{\ell(r^{-2})}.$$
\end{lemma}
\pf
The conclusion of this Lemma follows easily from
Theorem \ref{t:abofjfatorigin} above and Lemma 3.2 of \cite{SV}.
\qed

Using the L\'{e}vy system for $X$, we know that for every
bounded open subset $D$ and  every $f \ge 0$ and $x \in D$,
\begin{equation}\label{newls}
\E_x\left[f(X_{\tau_D});\,X_{\tau_D-} \not= X_{\tau_D}  \right]
=  \int_{\overline{D}^c} \int_{D}
G_D(x,z) J(z-y) dz f(y)dy.
\end{equation}
For notational convenience, we define
\begin{equation}\label{PK}
K_D(x,y)\,:=  \int_{D}
G_D(x,z) J(z-y) dz, \qquad (x,y) \in D \times
\overline{D}^c.
\end{equation}
Thus \eqref{newls} can be simply written as
$$
\E_x\left[f(X_{\tau_D});\,X_{\tau_D-} \not= X_{\tau_D}  \right]
=\int_{\overline{D}^c} K_D(x,y)f(y)dy.
$$
Using the continuities of $G_D$ and $J$, one can easily check
that $K_D$ is continuous on $D \times
\overline{D}^c$.

As a consequence of Lemma
\ref{l:tau}-\ref{l:tau2} and \eqref{PK}, we get the
following proposition.

\begin{prop}\label{p:Poisson1}
There exist  $C_1, C_2>0$ such that for every $r \in (0, r_2)$ and $x_0
\in \R^d$,
\begin{equation}\label{P1}
 K_{B(x_0,r)}(x,y) \,\le\, C_1 \, j(|y-x_0|-r) \frac{r^{\alpha/2}}
{(\ell(r^{-2}))^{1/2}}\frac{(r-|x-x_0|)^{\alpha/2}}
{(\ell((r-|x-x_0|)^{-2}))^{1/2}},
\end{equation}
for all $(x,y) \in B(x_0,r)\times \overline{B(x_0,r)}^c$ and
\begin{equation}\label{P2}
K_{B(x_0,r)}(x_0,y) \,\ge\, C_2\, J(2(y-x_0))\frac{r^\alpha}{\ell(r^{-2})},
\qquad\forall
y \in \overline{B(x_0,r)}^c.
\end{equation}
\end{prop}
\pf
Without loss of generality, we assume $x_0=0$.
For $z \in B(0, r)$ and $y \in \overline{B(0,r)}^c$,
$$|y|-r \le |y|-|z| \le |z-y| \le |z|+|y| \le r +|y| \le 2|y|.
$$
Thus by the monotonicity of $J$,
$$
J( 2y) \,\le\, J(z-y) \, \le \, j(|y|-r).     \qquad (z,y) \in B(0,r) \times
\overline{B(0,r)}^c.$$
Applying the above inequality and Lemma
\ref{l:tau}-\ref{l:tau2} to \eqref{PK}, we have proved the proposition.
\qed

\begin{prop}\label{p:Poisson2}
For every $a \in (0,1)$, there exists $C=C(a)>0$ such that for every
$r \in (0, r_2)$, $x_0 \in \R^d$ and $x_1, x_2 \in B(x_0, ar)$,
$$
K_{B(x_0,r)}(x_1,y) \,\le\, C  K_{B(x_0,r)}(x_2,y), \qquad
y \in  \overline{B(x_0,r)}^c.
$$
\end{prop}

\pf
This follows easily from the Harnack inequality (Theorem \ref{T:Har})
and the continuity of $K_{B(x_0, r)}$.
For details, see the proof of Lemma 4.2 in \cite{Sz2}.\qed

As an immediate consequence of Theorem \ref{t:abofjfatorigin}, we have

\begin{lemma}\label{l:J}
There exists $r_3 \in (0,r_2]$ such that for every
$y\in \R^d$ with $|y|\le r_3$,
$$
\frac{\alpha\Gamma((d+\alpha)/2)}{2^{2-\alpha}\pi^{d/2}\Gamma(1-
\alpha/2)}
\frac{\ell(|y|^{-2})}{|y|^{d+\alpha}} \,\le\, J(y) \,\le\,
\frac{2^{\alpha}\alpha\Gamma((d+\alpha)/2)}{\pi^{d/2}\Gamma(1-
\alpha/2)}\frac{\ell(|y|^{-2})}{|y|^{d+\alpha}}.
$$
\end{lemma}

The next inequalities will be used several times
in the remainder of this paper.

\begin{lemma}\label{l:l}
There exist $r_4 
\in (0,r_3]$ and $C
>0$
such that
\begin{equation}\label{el1}
\frac{s^{\alpha/2}}{\left(\ell(s^{-2})\right)^{1/2}} \,\le \, C \,
\frac{r^{\alpha/2}}{\left(\ell(r^{-2})\right)^{1/2}}, \qquad \forall\,
0<s<r\le 4 r_4,
\end{equation}
\begin{equation}\label{el2}
\frac{s^{1-\alpha/2}}{\left(\ell(s^{-2})\right)^{1/2}} \,\le \, C \,
\frac{r^{1-\alpha/2}}{\left(\ell(r^{-2})\right)^{1/2}}, \qquad  \forall\,
0<s<r\le 4r_4,
\end{equation}
\begin{equation}\label{el7}
s^{1-\alpha/2} \,{\left(\ell(s^{-2})\right)^{1/2}} \,\le \, C \,
r^{1-\alpha/2}\,{\left(\ell(r^{-2})\right)^{1/2}}, \qquad  \forall\,
0<s<r\le 4r_4,
\end{equation}
\begin{equation}\label{el3}
\int^{\infty}_r \frac{\left(\ell(s^{-2})\right)^{1/2}}{s^{1+\alpha/2}}ds
 \,\le \, C \,
\frac{\left(\ell(r^{-2})\right)^{1/2}}{r^{\alpha/2}}, \qquad \forall\,
0<r\le 4 r_4,
\end{equation}
\begin{equation}\label{el6}
\int^{r}_0 \frac{\left(\ell(s^{-2})\right)^{1/2}}{s^{\alpha/2}}ds
 \,\le \, C \,
\frac{\left(\ell(r^{-2})\right)^{1/2}}{r^{\alpha/2-1}}, \qquad \forall\,
0<r\le 4 r_4,
\end{equation}

\begin{equation}\label{el4}
\int^{\infty}_r \frac{\ell(s^{-2})}{s^{1+\alpha}}ds
 \,\le \, C \,
\frac{\ell(r^{-2})}{r^{\alpha}}, \qquad \forall\,  0<r\le 4 r_4,
\end{equation}
\begin{equation}\label{el5}
\int_{0}^r \frac{\ell(s^{-2})}{s^{\alpha-1}}ds
 \,\le \, C \,
\frac{\ell(r^{-2})}{r^{\alpha-2}}, \qquad \forall\,  0<r\le 4 r_4
\end{equation}
and
\begin{equation}\label{el8}
\int_{0}^r \frac{s^{\alpha-1}}{\ell(s^{-2})}ds
 \,\le \, C \,
\frac{r^{\alpha}}{\ell(r^{-2})}, \qquad \forall\,  0<r\le 4 r_4.
\end{equation}
\end{lemma}
\pf
The first three inequalities follow easily from Theorem 1.5.3 of \cite{BGT},
while the last five from the 0-version of Theorem 1.5.11 of \cite{BGT}.
\qed

\begin{prop}\label{p:Poisson3}
For every $a \in (0,1)$, there exists $C=C(a)>0$ such that for every
$r \in (0, r_4]$ and $x_0 \in \R^d$,
$$ K_{B(x_0,r)}(x,y) \,\le\, C\,\frac{r^{\alpha/2-d}}{(\ell(r^{-2}))^{1/2}}
\frac{(\ell((|y-x_0|-r)^{-2}))^{1/2}}
{( |y-x_0|-r)^{\alpha/2}}\qquad \forall
x\in  B(x_0, ar),\,  y \in \{r<|x_0-y| \le 2r\}.$$
\end{prop}

\pf
By Proposition \ref{p:Poisson2},
$$
K_{B(x_0,r)}(x,y) \le \frac{c_1}{r^d} \int_{B(x_0, a r)} K_{B(x_0,r)}(w,y) dw
$$
for some constant $c_1=c_1(a)>0$.
Thus from Lemma \ref{l:tau} and (\ref{P1})  
we have that
\begin{eqnarray*}
K_{B(x_0,r)}(x,y)&\le& \frac{c_2}{r^d}\int_{B(x_0, r)}\int_{B(x_0, r)}
G_{B(x_0,r)}(w,z)J(z-y) dz dw \\
&=& \frac{c_2}{r^d}\int_{B(x_0, r)}
\E_z[\tau_{B(x_0,r)}]J(z-y) dz\\
& \le& \frac{c_3}{r^d}  \frac{r^{\alpha/2}}{(\ell(r^{-2}))^{1/2}}
\int_{B(x_0, r)}\frac{(r-|z-x_0|)^{\alpha/2}}
{(\ell((r-|z-x_0|)^{-2}))^{1/2}} J(z-y)dz
\end{eqnarray*}
for some constants $c_2=c_2(a)>0$ and  $c_3=c_3(a)>0$.
Now applying Lemma \ref{l:J}, we get
$$
K_{B(x_0,r)}(x,y)
\, \le\,  \frac{c_4r^{\alpha/2-d}}{(\ell(r^{-2}))^{1/2}}
\int_{B(x_0, r)}\frac{(r-|z-x_0|)^{\alpha/2}}
{(\ell((r-|z-x_0|)^{-2}))^{1/2}} \frac{\ell(|z-y|^{-2})}
{|z-y|^{d+\alpha}}dz
$$
for some constant $c_4=c_4(a)>0$.
Since $r-|z-x_0| \le |y-z| \le 3r \le 3r_4$,  from \eqref{el1} we see that
$$
\frac{(r-|z-x_0|)^{\alpha/2}}
{(\ell((r-|z-x_0|)^{-2}))^{1/2}} \, \le\, c_5
 \frac{(|y-z|)^{\alpha/2}}
{(\ell(|y-z|^{-2}))^{1/2}}
$$
for some constant $c_5>0$.
Thus we have
\begin{eqnarray*}
K_{B(x_0,r)}(x,y)
& \le&   \frac{c_6r^{\alpha/2-d}}{(\ell(r^{-2}))^{1/2}}
\int_{B(x_0, r)}
\frac{(\ell(|z-y|^{-2}))^{1/2}}{|z-y|^{d+\alpha/2}}dz\\
& \le&   \frac{c_6r^{\alpha/2-d}}{(\ell(r^{-2}))^{1/2}}
\int_{B(y, |y-x_0|-r)^c}
\frac{(\ell(|z-y|^{-2}))^{1/2}}{|z-y|^{d+\alpha/2}}dz\\
&\le &
\frac{c_7r^{\alpha/2-d}}{(\ell(r^{-2}))^{1/2}}
\int_{|y-x_0|-r}^{\infty} \frac{\left(\ell(s^{-2})\right)^{1/2}}
{s^{1+\alpha/2}}ds
\end{eqnarray*}
for some constants $c_6=c_6(a)>0$ and  $c_7=c_7(a)>0$.
Using \eqref{el3} in the above equation, we conclude that
$$
K_{B(x_0,r)}(x,y)
\,\le
\,  \frac{c_8r^{\alpha/2-d}}{(\ell(r^{-2}))^{1/2}}
\frac{(\ell((|y-x_0|-r)^{-2}))^{1/2}}
{( |y-x_0|-r)^{\alpha/2}}
$$
for some constant $c_8=c_8(a)>0$.
\qed

\section{Boundary Harnack Principle}
In this section, we give the proof of the boundary Harnack principle for $X$.

Using an argument similar to the first part of the proof of Lemma 3.3 in \cite{SW}
and using Lemma \ref{l:J} and \eqref{el4}-\eqref{el5} above
we can easily get the following lemma. We skip the details.

\begin{lemma}\label{l2.1}
There exists $C>0$ such that for any
$r\in (0, r_4)$ and any open set $D$ with $D\subset B(0,
r)$ we have
$$
\P_x\left(X_{\tau_D} \in B(0, r)^c\right) \,\le\, C\,r^{-\alpha}\,
\ell(r^{-2})\int_D
G_D(x,y)dy, \qquad x \in D\cap B(0, r/2).
$$
\end{lemma}

\begin{lemma}\label{l2.1_1}
There exists $C>0$ such that for any open set
$D$ with  $B(A, \kappa r)\subset D\subset
B(0, r)$ for some $r\in (0, r_4)$ and $\kappa\in (0, 1)$, we have
that for every $x \in D \setminus B(A, \kappa r)$,
$$
\int_{D} G_D(x,y)   dy  \, \le\,C\, r^{\alpha} \,\kappa^{-d-\alpha/2}\,
\frac1{\ell((4r)^{-2})}\left(1+
\frac{\ell((\frac{\kappa r}{2})^{-2})}{\ell((4r)^{-2})}\right)
\P_x\left(X_{\tau_{D\setminus B(A, \kappa r)}} \in B(A, \kappa
r)\right).
$$
\end{lemma}

\pf Fix a point $x\in D\setminus B(A, \kappa r)$ and let $B:=B(A,
\frac{\kappa r}2)$. Note that, by the harmonicity of $G_D(x,\,\cdot\,)$ in
$D\setminus \{x\}$ with respect to $X$, we have
\[
G_D(x,A)\,\ge\,\int_{D\cap \overline{B}^c}K_B(A, y)G_D(x,y)dy
\,\ge\,\int_{D\cap B(A, \frac{3\kappa r}4)^c}K_B(A, y)G_D(x,y)dy.
\]

Since $\frac{3\kappa r}4\le |y-A|\le 2r$ for $y\in B(A,
\frac{3\kappa r}4)^c\cap D$ and $j$ is a decreasing function,
it follows from \eqref{P2} in Proposition
\ref{p:Poisson1} and Lemma \ref{l:J} that
\begin{eqnarray*}
G_D(x,A) &\ge& c_1\, \frac{(\frac{\kappa r}{2})^\alpha}{\ell\left((
\frac{\kappa r}{2})^{-2}\right)}\int_{D \cap B(A, \frac{3\kappa
r}4)^c}G_D(x,y)J(2(y-A)) dy\\
&\ge& c_1\, j(4r)\, \frac{(\frac{\kappa r}{2})^\alpha}{\ell\left((
\frac{\kappa r}{2})^{-2}\right)}\int_{D \cap B(A, \frac{3\kappa
r}4)^c}G_D(x,y) dy\\
&\ge& c_2\, \kappa^\alpha \,r^{-d}\, \frac{\ell((4r)^{-2})}{\ell((
\frac{\kappa r}{2})^{-2})}\int_{D \cap B(A, \frac{3\kappa
r}4)^c}G_D(x,y) dy,
\end{eqnarray*}
for some positive constants $c_1$ and $c_2$. On the other hand, applying Theorem
\ref{T:Har} we get
\[
\int_{B(A, \frac{3\kappa r}4)} G_D(x,y)   dy\le c_3 \int_{B(A,
\frac{3\kappa r}4)}  G_D(x,A)dy \,\le\,c_4\,r^{d}\,\kappa^d G_D(x,A),
\]
for some positive constants $c_3$ and $c_4$.
Combining these two estimates we get that
\begin{equation}\label{efe1}
\int_{D} G_D(x,y)   dy    \,\le\, c_5\,\left(r^{d}\kappa^d+r^{d}
\kappa^{-\alpha}\frac{\ell((\frac{\kappa r}{2})^{-2})}
{\ell((4r)^{-2})}\right)\,
G_D(x,A)
\end{equation}
for some constant $c_5>0$.

Let $\Omega=D\setminus B(A, \frac{\kappa r}2)$. Note that for any
$z\in B(A, \frac{\kappa r}4)$ and $y\in \Omega$, $ 2^{-1}|y-z|\le
|y-A|\le 2|y-z|. $ Thus we get from (\ref{PK}) that for $z\in B(A,
\frac{\kappa r}4)$,
\begin{equation}\label{e:KK1}
c_6^{-1}K_{\Omega}(x, A)  \,\le \,K_{\Omega}(x, z) \,\le\,
c_6K_{\Omega}(x, A)
\end{equation}
for some $c_6>1$.
Using the harmonicity of $G_D(\cdot, A)$ in $D\setminus\{A\}$ with
respect to $X$, we can split $G_D(\cdot, A)$ into two parts:
\begin{eqnarray*}
&&G_D(x, A)
=\E_x \left[G_D(X_{\tau_{\Omega}},A)\right]\\
&&=\E_x \left[G_D(X_{\tau_{\Omega}},A):\,X_{\tau_{\Omega}} \in B(A,
\frac{\kappa r}4)  \right]\,+\, \E_y
\left[G_D(X_{\tau_{\Omega}},A):\,X_{\tau_{\Omega}}
\in \{\frac{\kappa r}4\le |y-A|\le \frac{\kappa r}2\}\right]\\
&& :=I_1+I_2.
\end{eqnarray*}
Using (\ref{e:KK1})  and \eqref{G2}, we have
$$
I_1 \,\le\, c_6\,K_{\Omega}(x,A) \int_{B(A, \frac{\kappa
r}4)}G_D(y, A)dy \,\le\, c_7 \,K_{\Omega}(x,A) \int_{B(A,
\frac{\kappa r}4)}\frac{1}{|y-A|^{d-\alpha}} \frac{dy}{\ell(|y-A|^{-2})}.
$$
for some constant $c_7>0$.
Since $|y-A|\le 4r \le 4 r_4$, by \eqref{el1},
\begin{equation}\label{efe}
\frac{|y-A|^{\alpha/2}}{\ell(|y-A|^{-2})} \,\le\, c_8 \,
\frac{(4r)^{\alpha/2}}{\ell((4r)^{-2})}
\end{equation}
for some constant $c_8>0$.
Thus
$$
I_1 \,\le\, c_7\, c_8\,K_{\Omega}(x,A) \int_{B(A,
\frac{\kappa r}4)}\frac{1}{|y-A|^{d-\alpha/2}}
\frac{(4r)^{\alpha/2}}{\ell((4r)^{-2})}dy
\le\, c_9\kappa^{\alpha/2}r^{\alpha}\frac1{\ell((4r)^{-2})}K_{\Omega}(x, A)
$$
for some constant $c_9>0$.
Now   using (\ref{e:KK1}) again, we  get
 $$
  I_1      \,\le\,
c_{10}\kappa^{\alpha/2-d}r^{\alpha-d}\frac1{\ell((4r)^{-2})}\int_{B(A,
\frac{\kappa r}4)} K_{\Omega}(x,
z)dz,
$$
for some constant $c_{10}>0$. On the other
hand, by \eqref{G2},
\begin{eqnarray*}
&&I_2 \,=\,
\int_{\{\frac{\kappa r}4\le |y-A|\le \frac{\kappa r}2\}}
G_{D}(y,A) \P_x(X_{\tau_{\Omega}} \in dy)   \\
&&\le
 c_{11}\int_{\{\frac{\kappa r}4\le |y-A|\le \frac{\kappa r}2\}}
\frac1{|y-A|^{d-\alpha}} \,\frac{1}{\ell(|y-A|^{-2})}
\P_x(X_{\tau_{\Omega}} \in dy)
\end{eqnarray*}
for some constant $c_{11}>0$.
Using \eqref{efe}, the above is less than or equal to
$$ c_{12}
\kappa^{\alpha/2-d}\,r^{\alpha-d} \, \frac1{\ell((4r)^{-2})}\P_x
\left(X_{\tau_{\Omega}} \in
\{\frac{\kappa r}4\le |y-A|\le \frac{\kappa r}2\}\right),
$$
for some constant $c_{12}>0$.
Therefore
$$
G_D(x, A) \,\le\, c_{13}\,
\kappa^{\alpha/2-d}\,r^{\alpha-d}\frac1{\ell((4r)^{-2})}\,
\P_x\left(X_{\tau_{\Omega}} \in B(A,
\frac{\kappa r}2)\right).
$$
for some constant $c_{13}>0$.
Combining the above with \eqref{efe1},   we get
\[ \int_{D} G_D(x,y)   dy
\,\le\,c_{14}\, r^{\alpha} \,\kappa^{-d-\alpha/2}\,
\frac1{\ell((4r)^{-2})}\left(1+
\frac{\ell((\frac{\kappa r}{2})^{-2})}{\ell((4r)^{-2})}\right)\P_x
\left(X_{\tau_{D\setminus
B(A, \frac{\kappa r}2)}} \in B(A, \frac{\kappa r}2)\right),
\]
for some constant $c_{14}>0$. It follows
immediately that
$$
\int_{D} G_D(x,y)   dy  \, \le\,c_{14}\, r^{\alpha} \,
\kappa^{-d-\alpha/2}\,\frac1{l((4r)^{-2})}\left(1+
\frac{\ell((\frac{\kappa r}{2})^{-2})}
{\ell((4r)^{-2})} \right)
\P_x\left(X_{\tau_{D\setminus B(A, \kappa r)}} \in B(A, \kappa
r)\right).
$$
\qed

Combining Lemmas \ref{l2.1}-\ref{l2.1_1} and using
the translation invariant property, we have
the following

\begin{lemma}\label{l2.3}
There exists $c_1>0$ such that for any open set
$D$ with $B(A, \kappa r)\subset D\subset
B(Q, r)$ for some $r\in(0, r_4)$ and $\kappa\in (0, 1)$,
we have that for every $ x\in D\cap B(Q,
\frac{r}2)$,
$$
\P_x\left(X_{\tau_{D}} \in B(Q, r)^c\right) \,\le\,
c_1\,\kappa^{-d-\alpha/2 }\, \frac{\ell(r^{-2})}{ \ell((4r)^{-2})}\,
\left(1+\frac{\ell((\frac{\kappa r}{2})^{-2})}{\ell((4r)^{-2})}
\right)  \P_x\left(X_{\tau_{D\setminus B(A, \kappa r) }} \in
B(A, \kappa r) \right).
$$
\end{lemma}

Let  $A(x, a,b):=\{ y \in \R^d: a \le |y-x| <b  \}.$

\begin{lemma}\label{l2.U}
Let $D$ be an open set and $0<2r<r_4$.
For every $Q \in \R^d$ and any positive function $u$ vanishing
 on  $D^c \cap B(Q, \frac{11}6)$,
there is a $\sigma\in (\frac{10}{6}r, \frac{11}{6}r)$
such that for any  $M \in (1,\infty)$ and $ x \in D
\cap B(Q, \frac{3}{2}  r)$,
$$
\E_x\left[u(X_{\tau_{D \cap B(Q, \sigma)}});
X_{\tau_{D \cap B(Q, \sigma)}} \in
A(Q, \sigma, M) \right]\\
\,\le\,
C\,\frac{r^{\alpha}}{\ell((2r)^{-2})}   \int_{A(Q, \frac{10r}6, M)}
J(y)u(y)dy
$$
for some constant $C=C(M)>0$ independent of $Q$ and $u$.
\end{lemma}

\pf Without loss of generality, we may assume that $Q=0$.
Note that by \eqref{el6}
\begin{eqnarray*}
&&\int^{\frac{11}{6}r}_{\frac{10}{6}r}\int_{A(0, \sigma, 2r)}
\ell((|y|- \sigma)^{-2})^{1/2} (|y|-\sigma)^{-{\alpha}/2}
u(y)dy d\sigma\\
 &&=\int_{A(0, \frac{10}{6}r , 2r)}
\int^{ |y| \wedge \frac{11}{6}r}_{\frac{10}{6}r}\ell((|y|- \sigma)^{-2})^{1/2}
(|y|-\sigma)^{-{\alpha}/2}d\sigma u(y)dy \\
&& \le c_1 \int_{A(0, \frac{10}{6}r , 2r)}
\left(\int^{ |y| - \frac{10}{6}r}_{0}\ell(s^{-2})^{1/2}
s^{-{\alpha}/2}ds \right)u(y)dy \\
&&\le c_2
\int_{A(0, \frac{10r}6, 2r)} \ell((|y|-  \frac{10r}6 )^{-2})^{1/2}
(|y|- \frac{10r}6)^{1-{\alpha}/2} u(y)dy
\end{eqnarray*}
for some positive constants $c_1$ and $c_2$.
Using \eqref{el7}, we get that there is a constant $c_3>0$ such that
$$
\int_{A(0, \frac{10r}6, 2r)}\ell((|y|-  \frac{10r}6 )^{-2})^{1/2}
(|y|- \frac{10r}6)^{1-{\alpha}/2} u(y)dy
\le c_3
\int_{A(0, \frac{10r}6, 2r)} \ell(|y|^{-2})^{1/2} |y|^{1-{\alpha}/2} u(y)dy,
$$
which is less than or equal to
$$
 c_4 \frac{r^{1-\alpha/2}}{\ell((2r)^{-2})^{1/2}}  \int_{A(0,
\frac{10r}6, 2r)} \ell(|y|^{-2}) u(y)dy
$$
for some constant $c_4>0$ by \eqref{el2}.
Thus, by taking $c_5>6 c_2 c_4$, we can conclude that there is a
$\sigma\in (\frac{10}{6}r, \frac{11}{6}r)$ such that
\begin{equation}\label{e:int}
\int_{A(0, \sigma, 2r)}\ell((|y|- \sigma)^{-2})^{1/2}\,
(|y|-\sigma)^{-{\alpha}/2}u(y)dy
\,\le\, c_5\,\frac{r^{-\alpha/2}}{\ell((2r)^{-2})^{1/2}}
\int_{A(0, \frac{10r}6, 2r)} \ell(|y|^{-2}) u(y)dy.
\end{equation}

Let  $x \in D \cap B(0,  \frac{3}{2}  r)$.
Note that, since $X$ satisfies the hypothesis ${\bf H}$ in \cite{Sz1},
by Theorem 1 in \cite{Sz1} we have
\begin{eqnarray*}
&& \E_x\left[u(X_{\tau_{D \cap B(0, \sigma)}}); X_{\tau_{D \cap
B(0, \sigma)}} \in
A(0, \sigma, M) \right]\\
&&= \E_x\left[u(X_{\tau_{D \cap B(0, \sigma)}}); X_{\tau_{D \cap
B(0, \sigma)}} \in
A(0, \sigma, M), \, \tau_{D \cap
B(0, \sigma)} =\tau_{B(0, \sigma)}  \right]\\
&&= \E_x\left[u(X_{\tau_{ B(0, \sigma)}}); X_{ \tau_{B(0, \sigma)}} \in
A(0, \sigma, M), \, \tau_{D \cap B(0, \sigma)} =
\tau_{B(0, \sigma)}  \right]\\
&&\le \E_x\left[u(X_{\tau_{ B(0, \sigma)}}); X_{\tau_{B(0, \sigma)}} \in
A(0, \sigma, M)  \right]
\,=\,\int_{A(0, \sigma, M)}K_{B(0, \sigma)}(x, y)u(y)dy.
\end{eqnarray*}
In the first equality above we have used the fact that
$u$ vanishes on $D^c \cap B(0, \sigma)$.
Since  $ \sigma  <2r < r_4$, from \eqref{P1} in Proposition \ref{p:Poisson1},
Proposition \ref{p:Poisson3} and Lemma \ref{l:J}
we have
\begin{eqnarray*}
&& \E_x\left[u(X_{\tau_{D \cap B(0, \sigma)}}); X_{\tau_{D \cap
B(0, \sigma)}} \in
A(0, \sigma, M) \right]
\,\le\,
 \int_{  A(0, \sigma, M)   }
K_{B(0, \sigma)}(x, y)u(y)dy\\
&&\le\,
c_6 \int_{A(0, \sigma, 2r)}
\frac{\sigma^{\alpha/2-d}}{\left(\ell(\sigma^{-2})\right)^{1/2}}
\frac{(\ell((|y|-\sigma)^{-2}))^{1/2}}
{( |y|-\sigma)^{\alpha/2}} u(y)dy\\&&~~ +
c_6 \int_{A(0, 2r , M )}j(|y|-\sigma) \frac{\sigma^{\alpha/2}}
{(\ell(\sigma^{-2}))^{1/2}}\frac{(\sigma-|x|)^{\alpha/2}}
{(\ell((\sigma-|x|)^{-2}))^{1/2}}
 u(y)dy
\end{eqnarray*}
for some constant $c_6>0$.
For $y \in A(0,  2r , M)   $,
$ \frac1{12}|y|   \le |y|-\sigma $
and $ \sigma-|x|\le\sigma \le {2r}$. Thus by \eqref{el1} and
the monotonicity of $j$, we have
$$
j(|y|-\sigma)
\frac{\sigma^{\alpha/2}}{(\ell(\sigma^{-2}))^{1/2}}
\frac{(\sigma-|x|)^{\alpha/2}}
{(\ell((\sigma-|x|)^{-2}))^{1/2}}
 \,\le\, c_7 j(\frac{|y|}{12})
\frac{r^{\alpha}}{\ell((2r)^{-2})}
$$
for some constant $c_7>0$.
Thus by \eqref{H:1}, we get
$$
j(|y|-\sigma)
\frac{\sigma^{\alpha/2}}{(\ell(\sigma^{-2}))^{1/2}}
\frac{(\sigma-|x|)^{\alpha/2}}
{(\ell((\sigma-|x|)^{-2}))^{1/2}}
 \,\le\, c_8    j(|y|)
\frac{r^{\alpha}}{\ell((2r)^{-2})}
$$
for some constant $c_8=c_8(M)>0$.
On the other hand, by \eqref{el1} and \eqref{e:int}, there exist
positive constants $c_9$ and $c_{10}$ such that
\begin{eqnarray*}
&&\int_{A(0, \sigma, 2r)}
\frac{\sigma^{\alpha/2-d}}{\left(\ell(\sigma^{-2})\right)^{1/2}}
\frac{(\ell((|y|-\sigma)^{-2}))^{1/2}}
{( |y|-\sigma)^{\alpha/2}} u(y)dy\\
&\le& (\frac{10r}{6})^{-d}\frac{\sigma^{\alpha/2}}
{\left(\ell(\sigma^{-2})\right)^{1/2}}
\int_{A(0, \sigma, 2r)}
   \frac{(\ell((|y|-\sigma)^{-2}))^{1/2}}
{( |y|-\sigma)^{\alpha/2}} u(y)dy\\
&\le & c_9 r^{-d}\frac{(2r)^{\alpha/2}}{\left(\ell((2r)^{-2})\right)^{1/2}}
\,\frac{r^{-\alpha/2}}{\left(\ell((2r)^{-2})\right)^{1/2}}
\int_{A(0, \frac{10r}6, 2r)} \ell(|y|^{-2}) u(y)dy \\
&\le & c_{10} \frac{r^{\alpha}}{\ell((2r)^{-2})}
\int_{A(0, \frac{10r}6, 2r)}\ell(|y|^{-2}) |y|^{-d-\alpha} u(y) dy,
\end{eqnarray*}
which is less than or equal to
$$
 c_{11} \frac{r^{\alpha}}{\ell((2r)^{-2})}
\int_{A(0, \frac{10r}6, 2r)}J(y) u(y) dy,
$$
for some constants $c_{11}>0$ by Lemma \ref{l:J}.
Hence
$$
\E_x\left[u(X_{\tau_{D \cap B(0, \sigma)}});
X_{\tau_{D \cap B(0, \sigma)}} \in
A(0, \sigma, M) \right]\\
\,\le\,
c_{12}\,\frac{r^{\alpha}}{\ell((2r)^{-2})} \int_{A(0, \frac{10r}6, M)}
J(y)u(y)dy
$$
for some constant $c_{12}=c_{12}(M)>0$.
\qed

\begin{lemma}\label{l2.2}
Let $D$ be an open set.
Assume that $B(A, \kappa r)\subset D\cap B(Q, r)$
for some $0<r<2r_4$ and $\kappa\in (0, \frac12]$.
Suppose that $u\ge0$ is regular harmonic in $D\cap
B(Q, 2r)$ with respect to $X$
and $u=0$ in $(D^c\cap B(Q, 2r)) \cup B(Q, M)^c$. If $w$ is
a regular harmonic function with respect to $X$ in $D\cap B(Q, r)$
such that
\[
w(x)=\left\{
\begin{array}{ll}
u(x), & x\in B(Q, \frac{3r}2)^c\cup (D^c\cap B(Q, r)),\\
0, & x \in A(Q, r, \frac{3r}2),
\end{array}\right.
\]
then
\[
u(A) \,\ge \,w(A)\, \ge\, C\,\kappa^{\alpha} \frac{\ell((2r)^{-2})}
{\ell((\kappa r)^{-2})} \,u(x), \quad x \in D
\cap B(Q,\frac32 r)
\]
for some constant $C=C(M)>0$.
\end{lemma}

\pf
Without loss of generality, we may assume $Q=0$
and  $x \in D \cap B(0,\frac32 r)$.
The left hand side inequality in the conclusion of the lemma
is obvious, so we only need to prove the right hand side
inequality.
Since $u$ is regular harmonic in $D\cap B(0, 2r)$ with respect
to $X$ and $u=0$ on
 $B(0, M)^c$, we know from Lemma \ref{l2.U} that
there exists $\sigma\in (\frac{10r}6, \frac{11r}6)$ such that
$$
u(x)=
\E_x\left[u(X_{\tau_{D \cap B(0, \sigma)}}); \,X_{\tau_{D \cap
B(0, \sigma)}} \in
A(0, \sigma, M) \right]
\le
c_1\frac{r^{\alpha}}{\ell((2r)^{-2})} \int_{A(0, \frac{10r}6, M)}
J(y)u(y)dy
$$
for some constant $c_1=c_1(M)>0$.
On the other hand,
by \eqref{P2} in Proposition \ref{p:Poisson1},  we have  that
\begin{eqnarray*}
&&w(A)\,=\,\int_{A(0, \frac{3r}2 , M )}K_{D\cap B(0, r)}(A, y)u(y)dy
\,\ge\, \int_{A(0, \frac{3r}2, M )}K_{B(A, \kappa r)}(A, y)u(y)dy\\
&&\ge\, c_2\int_{A(0, \frac{3r}2,  M)}    J(2(A-y))
\frac{(\kappa r)^\alpha}{\ell((\kappa r)^{-2})}  u(y)dy
\end{eqnarray*}
for some constant $c_2>0$.
Note that $|y-A|\le  2|y|$ in $ A(0, \frac{3r}2, M) $.
Hence by the monotonicity of $j$ and \eqref{H:1},
\[
w(A)\,\ge\, c_2\,\frac{(\kappa r)^\alpha}{\ell((\kappa r)^{-2})}
\int_{A(0, \frac{3r}2,  M)} j( 4|y|)
u(y) dy\,\ge\, c_3\,\frac{(\kappa r)^\alpha}{\ell((\kappa r)^{-2})}
\int_{A(0, \frac{3r}2,  M)} J( y)
u(y) dy
\]
for some constant $c_3=c_3(M)>0$.
Therefore
$$
w(A)\,\ge \,c_4 \,\kappa^{\alpha}  \frac{\ell((2r)^{-2})}
{\ell((\kappa r)^{-2})}   \,u(x)
$$
for some constant $c_4=c_4(M)>0$.
\qed

We recall the definition of
$\kappa$-fat set from \cite{SW}.

\begin{defn}\label{fat}
Let $\kappa \in (0,1/2]$. We say that an open set $D$ in
$\R^d$ is $\kappa$-fat if there exists $R>0$ such that for each
$Q \in \partial D$ and $r \in (0, R)$,
$D \cap B(Q,r)$ contains a ball $B(A_r(Q),\kappa r)$.
The pair $(R, \kappa)$ is called the characteristics of
the $\kappa$-fat open set $D$.
\end{defn}

Note that all Lipschitz domain and all non-tangentially accessible
domain (see \cite{JK} for the definition) are $\kappa$-fat.
Moreover, every {\it John domain} is  $\kappa$-fat (see Lemma 6.3 in
\cite{MV}). The boundary of a $\kappa$-fat open set can be highly
nonrectifiable and, in general, no regularity of its boundary can be
inferred.  Bounded $\kappa$-fat open set may be disconnected.

Since $l$ is slowly varying at $\infty$, we get the Carleson's estimate from Lemma \ref{l2.2}.

\begin{corollary}\label{c:Carl}
Suppose that $D$ is  a bounded   $\kappa$-fat open set with
the characteristics $(R, \kappa)$. There exists a constant
$R_1$ such that if $r \le R_1$,
$Q\in \partial D$, $u\ge0$ is regular harmonic in $D\cap
B(Q, 2r)$ with respect to $X$
and $u=0$ in $(D^c\cap B(Q, 2r)) \cup B(Q, M)^c$,
then
\[
u(A) \,\ge \,u(x), \quad x \in D
\cap B(Q,\frac32 r)
\]
for some constant $C=C(M)>0$.
\end{corollary}

The next theorem is a boundary Harnack principle for bounded
$\kappa$-fat open set
and it is the main result of this section. Maybe a word of
caution is in order here. The boundary Harnack principle here is a
little different from the ones proved in \cite{B} and \cite{SW} in
the sense that in the boundary Harnack principle below we require
our harmonic functions to
vanish on the complement of some large open ball, which is a little weaker than
assuming that it vanishes on the whole complement of the open set.
However, this will not affect our application
later since we are mainly interested in the case when the harmonic
functions are given by the Green functions.

\begin{thm}\label{BHP}
Suppose that $D$ is  a bounded   $\kappa$-fat open set with
the characteristics $(R, \kappa)$. There exists a constant
$r_5:=r_5(D, \alpha, l) \le r_4 \wedge R$ such that if $2r \le r_5$ and
$Q\in \partial D$,
then for any nonnegative functions $u, v$ in $\RR^d$
which are regular harmonic in $D\cap B(Q, 2r)$ with respect to $X$
and vanish 
in $D^c \cap B(Q, 2r) \cup B(Q, M)^c$, 
we have
$$
C^{-1}\,\frac{u(A_r(Q))}{v(A_r(Q))}\,\le\,
\frac{u(x)}{v(x)}\,\le C\,\frac{u(A_r(Q))}{v(A_r(Q))},
\qquad x\in D\cap B(Q, \frac{r}2),
$$
for some constant
$C=C(D,M)>1$.
\end{thm}

\pf
Since $l$ is slowly varying at $\infty$, there exists a
constant $r_5:=r_5(D, \alpha, l) \le r_4 \wedge R$ such that for every
$2r \le r_5$,
\begin{equation}\label{lll}
\max \left(\frac{\ell(r^{-2})}{ \ell((\kappa r)^{-2})    },\,
\frac{\ell((2r)^{-2})}{ \ell((4r)^{-2})},\,
\frac{\ell((\frac{\kappa r}{2})^{-2})}{\ell((4r)^{-2})},\,
\frac{\ell((\kappa r)^{-2})}{\ell((2r)^{-2})}
\right) \,\le\, 2.
\end{equation}

Fix $2r \le r_5$ throughout this proof.
Without loss of generality we may assume that $Q=0$ and $u(A_r(0))=v(A_r(0))$.
For simplicity, we will write $A_r(0)$ as $A$ in the remainder
of this proof.
Define $u_1$ and $u_2$ to be regular harmonic functions
in $D\cap B(0, r)$ with respect to $X$ such that
$$
u_1(x)=\left\{
\begin{array}{ll}
u(x), & r\le |x|<\frac{3r}2,\\
0, & x\in B(0, \frac{3r}2)^c\cup(D^c\cap B(0, r))
\end{array}
\right.
$$
and
$$
u_2(x)=\left\{
\begin{array}{ll}
u(x), & x\in B(0, \frac{3r}2)^c\cup(D^c\cap B(0, r)),\\
0, & r\le |x|<\frac{3r}2,
\end{array}
\right.
$$
and note that $u=u_1+u_2$. If $D\cap\{r\le |y|<\frac{3r}2\}$
is empty, then $u_1=0$ and the inequality (\ref{e2.6}) below
holds trivially.
So we assume $D\cap\{r\le |y|<\frac{3r}2\}$ is not empty.
Then by Lemma \ref{l2.2},
$$
u(y)\le c_1 \kappa^{-\alpha} \frac{\ell((\kappa r)^{-2})}{\ell((2r)^{-2})}
u(A), \qquad y\in D\cap B(0,
\frac{3r}2),
$$
for some constant $c_1=c_1(M)>0$.
For $x\in D\cap B(0, \frac{r}2)$, we have
\begin{eqnarray*}
u_1(x)&=& \E_x\left[u(X_{\tau_{D\cap B(0, r)}}): X_{\tau_{D\cap B(0, r)}}\in
D\cap \{r\le |y|<\frac{3r}2\}\right]\\
&\le&\left(\sup_{D\cap\{r\le |y|<\frac{3r}2\}}u(y)\right)
   \P_x\left( X_{\tau_{D\cap B(0, r)}}\in
D\cap \{r\le |y|<\frac{3r}2\}\right)              \\
&\le&\left(\sup_{D\cap\{r\le |y|<\frac{3r}2\}}u(y)\right)
   \P_x\left( X_{\tau_{D\cap B(0, r)}}\in
B(0,r)^c        \right)     \\
&\le&c_1\,\kappa^{-\alpha} \frac{\ell((\kappa r)^{-2})}{\ell((2r)^{-2})}
\,u(A) \,\P_x\left( X_{\tau_{D\cap B(0, r)}}\in
B(0,r)^c \right).
\end{eqnarray*}
Now using Lemma \ref{l2.3} and \eqref{lll} we have that for $ x\in D\cap
B(0, \frac{r}2)$,
\begin{eqnarray}
&&u_1(x)\\
&&\le\, c_2\,\kappa^{-d-\frac32\alpha }\, \frac
{\ell((\kappa r)^{-2})}{\ell((2r)^{-2})}\frac{\ell(r^{-2})}{\ell((4r)^{-2})}\,
\left(1+\frac{\ell((\frac{\kappa r}{2})^{-2})}{\ell((4r)^{-2})}
\right)\,u(A)\,\P_x\left(
X_{\tau_{(D\cap B(0,r))\setminus
B(A, \frac{\kappa r}2)}} \in B(A, \frac{\kappa r}2)\right) \nonumber\\
&&\le\,c_3 \,u(A)\,\P_x\left(
X_{\tau_{(D\cap B(0,r))\setminus
B(A, \frac{\kappa r}2)}} \in B(A, \frac{\kappa r}2)\right)
\label{e2.3}
\end{eqnarray}
for some positive constants $c_2=c_2(M)$ and $c_3=
c_3(M, \kappa)$.
Since $2r < r_4$, Theorem \ref{T:Har} implies that
\[
u(y)\,\ge\, c_4\,u(A), \qquad y\in B(A, \frac{\kappa r}2)
\]
for some constant $c_4>0$.
Therefore for $x\in D\cap
B(0, \frac{r}2)$
\begin{equation}\label{e2.4}
u(x) \,=\,  \E_x\left[u(X_{\tau_{(D\cap B(0, r))\setminus
B(A, \frac{\kappa r}2)}}) \right]     \,\ge\, c_4\,u(A)\,
\P_x\left(X_{\tau_{(D\cap B(0,r))\setminus
B(A, \frac{\kappa r}2)}} \in B(A, \frac{\kappa r}2)\right).
\end{equation}
Using (\ref{e2.3}), the analogue of (\ref{e2.4}) for $v$
and the assumption that $u(A)=v(A)$, we get that for $x\in D\cap
B(0, \frac{r}2)$,
\begin{equation}\label{e2.6}
u_1(x)\,\le \,c_3\,v(A)\,
\P_x\left(X_{\tau_{(D\cap B(0, r)) \setminus
B(A, \frac{\kappa r}2)}} \in B(A, \frac{\kappa r}2)\right)\,\le \,c_5\,v(x)
\end{equation}
for some constant $c_5=c_5(M, \kappa)>0.$
Since $u=0$ on $B(0,M)^c$, we have that
for $x\in D\cap B(0, r)$,
\begin{eqnarray*}
u_2(x)&=& \int_{A(0, \frac{3r}2,  M)}K_{D\cap B(0, r)}
(x, z)u(z)dz\\
&=&\int_{A(0, \frac{3r}2,  M )}
\int_{D\cap B(0, r)}
G_{D\cap B(0, r)}(x, y)  J(y-z)dy u(z)dz.
\end{eqnarray*}
Let
\[
s(x)\,:=\,\int_{D\cap B(0, r)}G_{D\cap B(0, r)}(x, y)dy.
\]
Note that for every $y \in B(0,r)$ and $z \in B(0, \frac{3r}2)^c$,
$$
\frac13|z| \,\le\, |z|-r \,\le\, |z|-|y| \, \le\, |y-z|\, \le \,
|y|+|z|\, \le \, r+|z| \le 2 |z|.$$
So by the monotonicity of $j$, for every $y \in B(0,r)$ and
$z \in B(0, \frac{3r}2)$,
$$  j(12|z|) \, \le \,j(2|z|) \, \le \,J(y-z) \, \le \,
j(\frac13|z|) \, \le \,j(\frac1{12}|z|).
$$
Using \eqref{H:1}, we have that, for every $y \in B(0,r)$ and
$z \in A(0, \frac{3r}2,  M )$,
$$  c_6^{-1} j(|z|) \, \le  \,J(y-z) \, \le \,c_6\,j(|z|)
$$
for some constant $c_6=c_6(M)>0$.
Thus we have
\begin{equation}\label{e2.5}
c_7^{-1}\,\le \,\frac{u_2(x)}{u_2(A)}/\frac{s(x)}{s(A)}\,\le \,c_7,
\end{equation}
for some constant $c_7=c_7(M)>1$.
Applying (\ref{e2.5}) to $u$ and $v$ and Lemma \ref{l2.2} to $v$ and
$v_2$, we obtain for
$x\in D\cap B(0, \frac{r}2)$,
\begin{equation}\label{e2.7}
u_2(x)\,\le\, c_7\,u_2(A)\,\frac{s(x)}{s(A)}\,\le\, c_{7}^2\,
\frac{u_2(A)}{v_2(A)}\,v_2(x)\,
\le\, c_{8}\, \kappa^{-\alpha} \frac{\ell((\kappa r)^{-2})}
{\ell((2r)^{-2})}\frac{u(A)}{v(A)}\,v_2(x)\,=\,
c_{8}\,\kappa^{-\alpha} \frac{\ell((\kappa r)^{-2})}{\ell((2r)^{-2})}\,v_2(x),
\end{equation}
for some constant $c_8=c_8(M)>0.$
Combining (\ref{e2.6}) and (\ref{e2.7}) and applying \eqref{lll}, we have
\[
u(x)\,\le\, c_{9} \,v(x), \qquad x\in D\cap B(0,
\frac{r}2),
\]
for some constant $c_{9}=c_{9}(M, \kappa)>0.$
\qed

\section{Martin Boundary and Martin Representation}\label{sec-mbr}

In this section we will always assume that $D$ is a bounded $\kappa$-fat
open set in $\R^d$ with the characteristics
$(R, \kappa)$. We are going to apply Theorem \ref{BHP} to study the
Martin boundary of $D$ with respect to $X$.

We recall from Definition \ref{fat} that  for each $Q \in \partial
D$ and $r \in (0, R)$, $A_r(Q)$ is a point in $D \cap B(Q,r)$
satisfying $B(A_r(Q),\kappa r)  \subset D \cap B(Q,r)$. From
Theorem \ref{BHP}, we get the following boundary Harnack
principle for the Green function of $X$ which will play an important
role in this section. Recall that $r_5 \le R$ is the constant defined in Theorem
\ref{BHP}.

\begin{thm}\label{BHP2}
There exists a constant $c=c(D, \alpha, \ell)
>1$ such that for any $Q \in \partial D$, $r \in (0,r_5)$
and $z,w \in D \setminus B(Q,2r)$, we have $$ c^{-1}\,
\frac{G_D(z, A_r(Q))}{G_D(w, A_r(Q))} \,\le\,
\frac{G_D(z,x)}{G_D(w, x)} \,\le\, c\, \frac{G_D(z,A_r(Q))}
{G_D(w, A_r(Q))} ,\quad x\in D\cap B\left(Q,\frac{r}{2}\right). $$
\end{thm}

Since $\ell$ is slowly varying at $\infty$, there exists a positive constant
$r_6:=r_6(\kappa, l) \le r_5$ such that for every $2r \le r_6$,
\begin{equation}\label{ll}
 \frac1{2}\,\le\, \min \left(\frac{ \ell((\frac{\kappa^2}{64}
r^{-2})}{\ell(r^{-2})},\,
\frac{ \ell((\frac{4}{\kappa^2} r^{-2})} {\ell(r^{-2})}
\right)\,\le \, \max\left(\frac{\ell((\frac{\kappa^2}{64} r^{-2})}{\ell(r^{-2})},\,
\frac{ \ell((\frac{4}{\kappa^2} r^{-2})} {l(r^{-2})}
\right) \,\le\, 2.
\end{equation}

\begin{lemma}\label{l:5B}
There exist positive constants $c=c(D, \alpha)$ and
$\gamma=\gamma(D, \alpha)< \alpha$ such that for any $Q\in \partial
D$ and $r\in (0, r_6)$, and nonnegative function $u$ which is
harmonic with respect to $X$ in $D \cap B(Q, r)$ we have
\begin{equation}\label{e:gamma}
u(A_r(Q))\,\le\, c\,\left(\frac2{\kappa}\right)^{\gamma k}\,
\frac{\ell\left( (\kappa/2)^{-2k}r^{-2}\right))}{\ell\left(r^{-2}\right))}
 u(A_{(\kappa/2)^{k}r}(Q)), \qquad k=0, 1, \dots .
\end{equation}
\end{lemma}

\pf
Without loss of generality, we may assume $Q=0$.
Fix $r <  r_6$ and
let
$$
 \eta_k\,:=\,\left(\frac{\kappa}2\right)^{k}r,
 \quad A_k\,:=\, A_{ \eta_k}(0)
\quad \mbox{ and } \quad B_k\,:=\,B(A_k,  \eta_{k+1}), \quad k=0,1, \dots.
$$
Note that the $B_k$'s are disjoint. So by the harmonicity of $u$, we have
$$
u(A_k)
\,\ge\, \sum_{l=0}^{k-1} \E_{A_k}\left[u(Y_{\tau_{B_k}}):\,
Y_{\tau_{B_k}} \in B_l \right]\\
\,=\, \sum_{l=0}^{k-1} \int_{B_l} K_{B_k}(A_k, z) u(z) dz.
$$
Theorem  \ref{T:Har} implies that
$$
 \int_{B_l} K_{B_k}(A_k, z) u(z) dz \,\ge\, c_0\, u(A_l)
\int_{B_l} K_{B_k}(A_k, z) dz
$$
for some constant $c_0=c_0(d, \alpha)>0$. Since dist$(A_k, B_l)
\le 2 \eta_{l}$, by \eqref{P2} in Proposition \ref{p:Poisson1} and
the monotonicity of $j$ we have
$$
K_{B_k}(A_k,z)\ge\, c_1\,  J(2(A_k-z)) \frac{(\eta_{k+1})^\alpha}
{\ell((\eta_{k+1})^{-2})}
\ge\, c_1\,  J(4\eta_l )  \frac{(\eta_{k+1})^\alpha}{\ell
((\eta_{k+1})^{-2})},
 \qquad
z \in B_l.
$$
Applying Lemma \ref{l:J} and \eqref{ll}, for every $z$ in $B_l$ we get
$$
K_{B_k}(A_k,z)\,\ge\,c_2 \,
\frac{(\eta_{k+1})^\alpha}{(4\eta_{l})^{d+\alpha}}
\frac{\ell((4\eta_{l})^{-2})}{\ell((\eta_{l+1})^{-2})}
\frac{\ell((\eta_{l+1})^{-2})}{\ell((\eta_{k+1})^{-2})}
\ge\, 2\, c_2\left(\frac{\kappa}{8}\right)^{d+\alpha}\frac{(\eta_{k+1})^\alpha}{(\eta_{l+1})^{d+\alpha}}
 \frac{\ell((\eta_{l+1})^{-2})}{\ell((\eta_{k+1})^{-2})}
$$
for some constant $c_2=
c_2(d, \alpha, \ell) > 0$. Thus we have
$$
\int_{B_l}K_{B_k}(A_k, z)dz \ge\, c_3\,  \frac{(\eta_{k+1})^\alpha}
{(\eta_{l+1})^{\alpha}}
\frac{\ell((\eta_{l+1})^{-2})}{\ell((\eta_{k+1})^{-2})}, \qquad z \in B_l
$$
for some constant $c_3=c_3(d, \alpha, \ell)>0$. Therefore,
$$
\left( \eta_k\right)^{-\alpha} u(A_k) \ell((\eta_{k+1})^{-2}) \,
\ge \, c_4 \sum_{l=0}^{k-1}
\left( \eta_l\right)^{-\alpha} u(A_l)\ell((\eta_{l+1})^{-2})
$$
for some constant $c_4=c_4(d, \alpha, \kappa, \ell)>0$. Let $a_k :=
 (\eta_k)^{-\alpha}u(A_k) \ell(\frac1{(\eta_{k+1})^2}) $ so that
$a_k \ge  c_4\sum_{l=0}^{k-1}  a_l$. By
induction, one can easily check that $ a_k  \ge c_5 (1+c_4/2)^{k} a_0$
for some constant $c_5=c_5(d, \alpha)>0$. Thus, with $ \gamma =
\alpha - {\ln(1+\frac{c_4}2)} (\ln (2/\kappa))^{-1}, $ we get
$$
u(A_r(Q))\,\le\, c\,\left(\frac2{\kappa}\right)^{\gamma k}\,
\frac{\ell\left( (\kappa/2)^{-2(k+1)}r^{-2}\right))}{\ell\left( (\kappa/2)^{-2}
r^{-2}\right))} u(A_{(\kappa/2)^{k}r}(Q)).
$$
Applying \eqref{ll}, we conclude that
(\ref{e:gamma}) is true. \qed

\begin{lemma}\label{l:la}
Suppose $Q \in \partial D$ and $r \in (0, r_5)$. If $w \in D\setminus
B(Q, r)$, then
\[
G_D(A_r(Q), w)\, \ge\, c\, \frac{\kappa^{\alpha}r^{\alpha}}
{\ell((\kappa r/2)^{-2})}\int_{B(Q, r)^c} J (\frac12(z- Q)) G_D(z, w)dz
\]
for some constant $c=c(D, \alpha, \ell)>0$.
\end{lemma}

\pf Without loss of
generality, we may assume $Q=0$. Fix $w \in D\setminus B(0, r)$ and
let $A:=A_r(0)$ and $u(\cdot) := G_D(\cdot, w)$. Since $u$ is
regular harmonic in $D\cap B(0, (1-\kappa/2)r)$ with respect to $X$,
we have
\begin{eqnarray*}
&&u(A) \,\ge\, \E_A \left[ u\left( X_{\tau_{D \cap
B(0,(1-\kappa/2)r)}}\right); X_{\tau_{D \cap B(0,(1-\kappa/2)r)}}
\in
B(0, r)^c \right]\\
&&= \int_{B(0, r)^c}K_{D\cap
B(0, (1-\kappa/2)r)}(A, z)u(z)dz\\
&&=\int_{B(0, r)^c         }
\int_{D\cap B(0,  (1-\kappa/2)r ) }
 G_{D\cap B(0,(1-\kappa/2)r )}(A,y)\,J(y-z) dyu(z)dz.
\end{eqnarray*}
Since $B(A, \kappa r/2)\subset D \cap B(0, (1-\kappa/2)r )$,
by the monotonicity of
the Green functions,
$$
G_{D\cap B(0,(1-\kappa/2)r )}(A,y) \, \ge \,
G_{B(A, \kappa r/2)}(A,y),
\quad y \in B(A, \kappa r/2).
$$
Thus
\begin{eqnarray*}
u(A) &\ge&  \int_{B(0, r)^c }
\int_{ B(A, \kappa r/2)  }
G_{B(A, \kappa r/2)}(A,y) J(y-z) dyu(z)dz\\
&=& \int_{B(0, r)^c }
K_{B(A, \kappa r/2)}(A, z)u(z)dz,
\end{eqnarray*}
which is greater than or equal to
$$
c_1\int_{B(0, r)^c}
J(2(z-A)) \frac{(\kappa r/2)^{\alpha}}{\ell((\kappa r/2)^{-2})}\,u(z)dz
$$
for some positive constant $c_1=c_1(d, \alpha, l)$ by \eqref{P2}
in Proposition \ref{p:Poisson1}. Note
that $|z-A|\le  2|z|$ for $z
\in B(0,r)^c $. Let $M:=$diam$(D)$. Hence
\begin{equation}\label{e:22}
u(A)\,\ge\, c_2\,\frac{\kappa^{\alpha}r^{\alpha}}
{\ell((\kappa r/2)^{-2})}\int_{A(0,r , M)}
{u(z)}  J(4 z) dz\,\ge\, c_3\,\frac{\kappa^{\alpha}r^{\alpha}}
{\ell((\kappa r/2)^{-2})}\int_{A(0,r , M)}
{u(z)}  J(\frac12 z) dz\end{equation}
for some constant $c_3=c_3(d, \alpha, \ell, M)>0$.
We have used \eqref{H:1} in the last inequality above. \qed

\begin{lemma}\label{l:14B}
There exist positive constants $c_1=c_1(D,\alpha,l)$  and
$c_2=c_2(D,\alpha, l)<1$ such that for any $Q \in \partial D $, $r\in
(0, r_6)$ and  $w \in D \setminus B(Q,2r/\kappa)$,  we have
$$
\E_x\left[G_D(X_{\tau_{D \cap B_k}}, w):\,X_{\tau_{D
\cap B_k}} \in B(Q, r)^c
 \right] \,\le\, c_1\,c_2^{k} \, G_D(x,w), \quad x \in D \cap B_k,
$$
where $B_k:=B(Q, (\kappa/2)^{k}r)$, $ k=0,1, \dots$.
\end{lemma}

\pf Without loss of generality, we may assume $Q=0$. Fix $r <r_6$
and  $w \in D \setminus B(0,4r)$.  Let $\eta_k:=(\kappa/2)^{k}r $,
$B_k:=B(0,\eta_k)$ and
$$
u_k(x)  \,:=\, \E_x\left[G_D(X_{\tau_{D \cap B_k}}, w);
X_{\tau_{D \cap B_k}}
 \in B(0, r)^c \right], \quad x \in D \cap B_k.
$$
Note that for $x\in D\cap B_{k+1}$
\begin{eqnarray*}
u_{k+1}(x) &=&  \E_x\left[G_D(X_{\tau_{D \cap B_{k+1}}}, w);\,
X_{\tau_{D \cap B_{k+1}}}
 \in B(0, r)^c \right] \\
&=&  \E_x\left[G_D(X_{\tau_{D \cap B_{k+1}}}, w);\, \tau_{D
\cap B_{k+1}} = \tau_{D \cap B_{k}}, ~ X_{\tau_{D \cap B_{k+1}}}
 \in B(0, r)^c \right] \\
&=&  \E_x\left[ G_D(X_{\tau_{D \cap B_k}}, w)  ;\, \tau_{D
\cap B_{k+1}} = \tau_{D\cap B_{k}}, ~ X_{\tau_{D \cap B_{k}}}
 \in B(0, r)^c  \right] \\
 &\le &
 \E_x\left[G_D(X_{\tau_{D \cap B_k}}, w);\, X_{\tau_{D \cap
B_{k}}}
 \in B(0, r)^c  \right].
\end{eqnarray*}
Thus
\begin{equation}\label{e:dec1}
u_{k+1}(x) \,\le\, u_{k}(x)\, ,\quad x\in D\cap B_{k+1}\, .
\end{equation}
Let
$A_k\,:=\,A_{\eta_k}(0) $ and $M:=$diam$(D)$.
Since $G_D(\,\cdot\,, w)$ is zero on $D^c$, we have
\begin{eqnarray*}
&&u_{k}(A_k) \,= \, \E_{A_k}\left[G_D(X_{\tau_{D \cap B_k}}, w);\,
X_{\tau_{D \cap B_{k}}}
 \in A(0, r, M)  \right] \\
&&\le  \E_{A_k}\left[G_D(X_{\tau_{ B_k}}, w)  ;\, X_{\tau_{
B_{k}}}
 \in A(0, r, M) \right]
\,\le\,  \int_{A(0, r, M)} K_{B_k}(A_k,z) G_D(z,w) dz.
\end{eqnarray*}
Since $r<r_4$, by  \eqref{P1} in Proposition \ref{p:Poisson1},
we get that for $z\in A(0, r, M)$,
$$
K_{B_k}(A_k,z) \,\,\le\, c_1 \,
\, J(|z|-\eta_k) \frac{\eta_k^{\alpha/2}}{(\ell(\eta_k^{-2}))^{1/2}}
\frac{(\eta_k-|A_k|)^{\alpha/2}}
{(\ell((\eta_k-|A_k|)^{-2}))^{1/2}}
$$
for some constant $c_1=c_1(D, \alpha)>0$ and $k=1,2, \dots$.
Since $ \eta_k-|A_k| \le\eta_k   \le r_6$, from \eqref{el1} we see that
$$
\frac{(\eta_k-|A_k|)^{\alpha/2}}
{(\ell((\eta_k-|A_k|)^{-2}))^{1/2}}   \,\le\, c \,
\frac{\eta_k^{\alpha/2}}{(\ell(\eta_k^{-2}))^{1/2}}.
$$
Thus
$$
K_{B_k}(A_k,z) \, \le\, c_2
\,J(|z|-\eta_k) \frac{\eta_k^{\alpha}}{\ell(\eta_k^{-2})}
$$
for some constant $c_2=c_2(D, \alpha, \ell)>0$ and $k=1,2, \dots$.
Therefore by the monotonicity of $j$
\begin{equation}\label{e:dec2}
 u_{k}(A_k) \,\le\, c_2\, \frac{\eta_k^{\alpha}}{\ell(\eta_k^{-2})}
\int_{A(0, r, M)}J(\frac12 z)
G_D(z,w)dz , \quad k=1,2, \dots.
\end{equation}
 From Lemma \ref{l:la}, we
have
\begin{equation}\label{e:dec3}
 G_D(A_0,w) \,\ge\, c_3 \frac{\kappa^{\alpha}r^{\alpha}}
{\ell((\kappa r/2)^{-2})}\int_{A(0, r, M)} J (\frac12z) G_D(z, w)dz
\end{equation}
for some constant $c_3=c_3(D, \alpha, \ell)>0$. Therefore (\ref{e:dec2}) and
(\ref{e:dec3}) imply that
$$
u_{k}(A_k) \,\le\, c_4\,\left(\frac{\kappa}{2}\right)^{k\alpha}
\frac{\ell\left((\kappa/2)^{-2} r^{-2}\right)}
{\ell\left((\kappa/2)^{-2k}r^{-2}\right)} G_D(A_0,w)
$$
for some constant $c_4=c_4(D, \alpha, \ell)>0$. On the other
hand, using Lemma \ref{l:5B}, we get
$$G_D(A_0,w) \,\le\,
c_5\,\,\left(\frac2{\kappa}\right)^{\gamma k}\,
\frac{\ell\left( (\kappa/2)^{-2k}r^{-2}\right))}{\ell\left( r^{-2}\right))}
    G_D(A_k,w)$$
for some constant $c_5=c_5(D,
\alpha)>0$. Thus by \eqref{ll}
$$ u_{k}(A_k) \,\le\,
c_6
\left(\frac2{\kappa}\right)^{-k(\alpha-\gamma)}G_D(A_k,w)
$$ for some constant $c_6=c_6(D, \alpha)>0$ and $k=1,2, \dots
$. By Theorem \ref{BHP2},
we have
$$
\frac{u_k(x)}{G_D(x, w)}    \,\le\,
\frac{u_{k-1}(x)}{G_D(x, w)}
\,\le\, c_6\, \frac{u_{k-1}(A_{k-1})}{G_D(A_{k-1}, w)}  \,\le\,
c_4c_5c_6\,\left(\frac2{\kappa}\right)^{-(k-1)(\alpha-\gamma)}
$$
for $k=1,2, \dots.$ \qed

Let $x_0\in D$ be fixed and set
\[
M_D(x, y):=\frac{G_D(x, y)}{G_D(x_0, y)}, \qquad x, y\in D,~
y\neq x_0.
\]
$M_D$ is called the Martin kernel of $D$ with respect to $X$.

Now the next theorem follows from Theorem
\ref{BHP2} and Lemma \ref{l:14B} (instead of  Lemma 13 and
Lemma 14 in \cite{B} respectively) in very much the same way as in
the case of symmetric stable processes in Lemma 16 of \cite{B} (with
Green functions instead of harmonic functions). We omit the details.

\begin{thm}\label{t2.2}
There exist positive constants $R_1$, $M_1$, $c$ and $\beta$
depending on $D$, $\alpha$ and $l$ such that for any $Q \in \partial D $,
$r < R_1$ and $z \in D \setminus B(Q, M_1 r)$, we have
\[
\left|M_D(z, x)-M_D(z, y)\right| \,\le\,
c\,\left(\frac{|x-y|}r\right)^{\beta}, \qquad x, y\in  D \cap B(Q,
r).
\]
In particular, the limit $\lim_{D \ni y\to w} M_D(x, y)  $ exists
for every $w\in \partial  D$.
\end{thm}

There is a compactification $D^M$ of $D$, unique up to a homeomorphism,
such that $M_D(x, y)$ has a continuous
extension to $D\times (D^M\setminus\{x_0\})$ and $M_D(\cdot, z_1)
=M_D(\cdot, z_2)$ if and only if $z_1=z_2$. (See, for instance,
 \cite{KW}.)
The set $\partial^MD=D^M\setminus D$ is called the Martin boundary
of $D$. For $z\in \partial^MD$, set $M_D(\cdot, z)$ to be zero in
$D^c$.

A positive harmonic function $u$ for  $X^D$ is minimal if, whenever
$v$ is a positive harmonic function for $X^D$ with $v\le u$ on $D$,
one must have $u=cv$ for some constant $c$. The set of points $z\in
\partial^MD$ such that $M_D(\cdot, z)$ is minimal harmonic
for $X^D$ is called the minimal Martin boundary of $D$.

For each fixed $z\in \partial D$ and $x\in D$,
let
\[
M_D(x, z):=\lim_{D\ni y\to z}M_D(x,y),
\]
which exists by Theorem \ref{t2.2}. For each $z\in
\partial D$, set $M_D(x, z)$
to be zero for $x\in D^c$.

\begin{lemma}\label{l:MH1}
For every $z \in \partial D$ and $B \subset \overline{B} \subset D$,
$M_D(X_{\tau_{B}} , z)$ is $\P_x$-integrable.
\end{lemma}
\pf
Take a sequence $\{z_m\}_{m \ge 1} \subset D\setminus
\overline{B}$ converging to $z$.
Since $M_D (\cdot , z_m )$ is regular harmonic for $X$ in $B$,
by Fatou's lemma and Theorem \ref{t2.2},
$$  \E_x \left[ M_D\left(X_{\tau_{B}} , z\right)\right]\,=\,
\E_x \left[ \lim_{m\to \infty} M_D\left(X_{\tau_{B}} ,
z_m\right) \right]\,\leq\, \liminf_{m\to \infty} M_D(x, z_m)\,=\,
M_D(x,z) \,<\,\infty .
$$
\qed

\begin{lemma}\label{l:MH2}
For every $z \in \partial D$ and  $x \in D$,
\begin{equation}\label{eqn:4.1}
M_D(x,z) \,=\, \E_x \left[M_D\left(X^D_{\tau_{B(x,r)}} ,
z\right)\right], \quad \mbox{ for
 every } 0<r<r_6 \wedge \frac12\rho_D(x).
\end{equation}
\end{lemma}

\pf Fix $z \in \partial D$, $x \in D$ and $r<r_6 \wedge
\frac12\rho_D(x)<R$.
let
$$
 \eta_m\,:=\,\left(\frac{\kappa}2\right)^{m}r\quad \mbox{ and } \
 \quad z_m\,:=\, A_{ \eta_m}(0), \quad m=0,1, \dots.
$$
Note that
$$ B(z_m , \, \eta_{m+1}) \,\subset \, B(z,\,  \frac12 \eta_m
) \cap D
 \,\subset  \,B(z, \,  \eta_m )\cap D  \,\subset  \,B(z, r)\cap D
 \,\subset \, D \setminus B(x,r)
$$
for all $m \ge0 $. Thus by the harmonicity of $M_D(\cdot , z_m )$,
we have
$$M_D(x,z_m) \,=\, \E_x \left[M_D\left(X_{\tau_{B(x,r)}} , z_m\right)\right].
$$

On the other hand, by Theorem \ref{BHP2}, there exist constants $m_0
\ge 0$ and $c_1>0$ such that for every $w \in D \setminus B (z,
\eta_m)$ and  $y \in D \cap B(z, \eta_{m+1})$,
$$
M_D(w, z_m) \,=\, \frac{G_D(w, z_m)}{G_D(x_0 , z_m )} \,\leq\,
c_1\, \frac{G_D(w, y)}{G_D(x_0 , y )} \,=\, c_1\, M_D(w,y),
\quad m\geq m_0.
$$
Letting $y \rightarrow z\in \partial D$ we get
\begin{equation}\label{eqn:4.2}
M_D(w, z_m) \,\leq\, c_1\, M_D (w, z), \quad m\geq m_0,
\end{equation}
for every $w \in D \setminus B(z, \eta_m).$

To prove (\ref{eqn:4.1}), it  suffices to show that
$\{ M_D(X_{\tau_{B(x,r)}}, z_m ) : m\geq m_0 \}$  is $\P_x $-uniformly
integrable.
Since $M_D(X_{\tau_{B(x,r)}}, z)$ is $\P_x$-integrable
by Lemma \ref{l:MH1}, for
any $\eps>0$, there is an $N_0>1$ such that
\begin{equation}\label{eqn:c4.3}
\E_x \left[M_D\left(X_{\tau_{B(x,r)}}, z\right) ;\,
M_D\left(X_{\tau_{B(x,r)}}, z\right) > N_0/c_1 \right] \,<\,
\frac{\eps}{4c_1}.
\end{equation}
Note that by (\ref{eqn:4.2}) and (\ref{eqn:c4.3})
\begin{eqnarray*}
&& \E_x \left[M_D\left(X_{\tau_{B(x,r)}}, z_m\right) ; \,
M_D\left(X_{\tau_{B(x,r)}}, z_m\right) > N_0 ~\mbox{ and }~
X_{\tau_{B(x,r)}} \in D \setminus B(z, \eta_m)\right] \\
&\leq & c_1\, \E_x \left[M_D\left(X_{\tau_{B(x,r)}}, z\right) ;
\, c_1 M_D\left(X_{\tau_{B(x,r)}}, z\right) > N_0 \right]
 \,<\,   c_1\,\frac{ \eps}{4c_1} \,=\, \frac{\eps}{4} .
\end{eqnarray*}
By \eqref{P1} in Proposition \ref{p:Poisson1}, we have for $m\geq
m_0$,
\begin{eqnarray*}
&&\E_x \left[M_D\left(X^D_{\tau_{B(x,r)}}, z_m\right) ; \,
X_{\tau_{B(x,r)}} \in D
 \cap B(z, \eta_m)\right]\,=\,\int_{D \cap B(z, \eta_m)} M_D(w, z_m)
K_{B(x,r)}(x,w) dw \\
&&\leq \, c_2\,\int_{D \cap B(z, \eta_m)} M_D(w, z_m)
j(|w-x|-r) \frac{r^{\alpha/2}}{(\ell(r^{-2}))^{1/2}}\frac{(r-|w|)^{\alpha/2}}
{(\ell((r-|w|)^{-2}))^{1/2}} dw
\end{eqnarray*}
for some $c_2=c_2(d, \alpha, l)>0$.
Since $|w-x| \ge |x-z|-|z-w|\ge \rho_D(x)- \eta_m \ge 2r -r=r$,
using the monotonicity of $J$ and \eqref{el1} to the above equation,
we see that
\begin{eqnarray}
&&\E_x \left[M_D\left(X^D_{\tau_{B(x,r)}}, z_m\right) ; \,
X_{\tau_{B(x,r)}} \in D
 \cap B(z, \eta_m)\right]\\&&\leq \, c_3\,j(r)
\frac{r^{\alpha}}{\ell(r^{-2})}\int_{D \cap B(z, \eta_m)} M_D(w, z_m)
 dw\nonumber\\
&&\leq \, c_4 \, \int_{B(z, \eta_m)} M_D(w, z_m) dw \,=\, c_4 \,
G_D(x_0, z_m)^{-1} \int_{B(z,  \eta_m)} G_D(w, z_m)
dw\label{e:new1}
\end{eqnarray}
for some $c_3=c_3(D, \alpha, \ell)>0$ and $c_4=c_4(D, \alpha, \ell, r)>0$.
Note that, by Lemma \ref{l:5B}, there exist $c_5=c_5(D, \alpha, \ell, m_0)
>0$, $c_6=c_6(D, \alpha, \ell, m_0, r)
>0$ and $\gamma < \alpha$ such that
\begin{eqnarray}
G_D(x_0, z_m)^{-1} &\le& c_5\, \left(\frac{\kappa}2\right)^{-\gamma m}\,
\frac{\ell\left( (\kappa/2)^{-2(m+1)} (\kappa/2)^{-2m_0}  r^{-2}\right))}
{\ell\left( (\kappa/2)^{-2} (\kappa/2)^{-2m_0}r^{-2}\right)}
G_D(x_0, z_{m_0})^{-1}\nonumber \\
 &\le&
c_6\,   \left(\frac{\kappa}2\right)^{-\gamma m}\,
\ell\left( (\kappa/2)^{-2m} (\kappa/2)^{-2(m_0+1)} r^{-2}\right)
\label{e:new2} .
\end{eqnarray}

On the other hand, by \eqref{G2}
\begin{eqnarray}
\int_{B(z, \eta_m)} G_D(w, z_m)
dw
&\le& c_7 \int_{B(z_m,  2\eta_m)} \frac{dw}{\ell(|w-z_m|^{-2})
|w-z_m|^{d-\alpha}} \nonumber\\
&\le& c_8  \int_{0}^{2\eta_m} \frac{s^{\alpha-1}}{\ell(s^{-2})}ds
 \,\le \, c_9 \,
\frac{ (\eta_m)^{\alpha}}{\ell((2\eta_m)^{-2})} \label{e:new3}
\end{eqnarray}
In the last inequality above, we have used \eqref{el8}.
It follows from (\ref{e:new1})-(\ref{e:new3})
that there exists $c_{10}=c_{10}(D, \alpha, \ell, m_0, r)
>0$  such that
$$
\E_x \left[M_D(X^D_{\tau_{B(x,r)}}, z_m);
X_{\tau_{B(x,r)}} \in D\cap B(z, 2r /m)\right]\le  c_{10}
\left(\frac{\kappa}2\right)^{(\alpha-\gamma )m}
\frac{\ell\left( (\kappa/2)^{-2m} (\kappa/2)^{-2(m_0+1)} r^{-2}\right)}
{\ell\left( (\kappa/2)^{-2m} ( 2r)^{-2}\right)}.$$
Since $\ell$ is slowly varying at $\infty$, we can take
$N=N(\eps, D, m_0, r)$ large enough so that  for $m\geq N$,
\begin{eqnarray*}
&& \E_x  \left[M_D\left(X_{\tau_{B(x,r)}}, z_m\right);\,
M_D\left(X_{\tau_{B(x,r)}}, z_m\right) > N\right]\\
&\leq&  \E_x  \left[M_D\left(X_{\tau_{B(x,r)}}, z_m\right) ; \,
X_{\tau_{B(x,r)}} \in D \cap B(z, 2r/m)\right]\\
&&+ \E_x  \left[M_D\left(X_{\tau_{B(x,r)}}, z_m\right) ; \,
M_D\left(X_{\tau_{B(x,r)}}, z_m\right) > N \,\mbox{ and }\,
X_{\tau_{B(x,r)}} \in D \setminus B(z, 2r/m)\right] \\
&<&   c_{10}\,
\left(\frac{\kappa}2\right)^{(\alpha-\gamma )m}\,
\frac{\ell\left( (\kappa/2)^{-2m} (\kappa/2)^{-2(m_0+1)}
r^{-2}\right)} {\ell\left( (\kappa/2)^{-2m} ( 2r)^{-2}\right)}\,+\,
\frac{\eps}{4}
\,\,<\,\,  \eps.
\end{eqnarray*}
As each $M_D(X_{\tau_{B(x,r)}}, z_m)$ is $\P_x$-integrable,
we conclude that
$\{ M_D(X_{\tau_{B(x,r)}}, z_m ) : m\geq m_0 \}$  is uniformly integrable
under $\P_x$. \qed

Using the fact that $\P_x(X_{\tau_{U}} \in \partial U)=0$ for every smooth open set $U$
(Theorem 1 in \cite{Sz1}), one can follow the proof of Theorem 2.2 of \cite{CS} or the proof
of Theorem 4.8 of \cite{KS2} and show that
The two lemmas above imply that $M_D(\cdot, z)$ is harmonic for $X$. We skip the details.

\begin{thm}\label{T:L4.3}
For every $z \in \partial D$,
the  function $x\mapsto M_D(\cdot, z)$ is harmonic in $D$
with respect to $X$.
\end{thm}

Recall that a point $z\in \partial D$ is said to be a regular
boundary point for $X$ if $\P_z(\tau_D=0)=1$ and an irregular
boundary point if $\P_z(\tau_D=0)=0$. It is well known that  if
$z\in \partial D$ is regular for $X$, then for any $x\in D$,
$G_D(x, y)\rightarrow 0$ as $y\rightarrow z$.

\begin{lemma}\label{lemma:boy0}
\begin{description}
\item{(1)} If
$z, w\in\partial D$, $z\neq w$ and $w$ is a regular boundary point for $Y$,
then $M_D(x, z)\to 0$ as
$x\to w$.\item{(2)}
The mapping
$(x, z)\mapsto M_D(x, z)$ is continuous on
$D\times\partial D$.
\end{description}
\end{lemma}

\pf Both of the assertions can be proved easily using our Theorems
\ref{BHP2} and \ref{t2.2}. We skip the proof since the
argument is almost identical to the one on  page 235 of \cite{B2}.
\qed

\begin{lemma}\label{l3.1}
Suppose that $h$ is
a bounded singular $\alpha$-harmonic function in
a bounded open set $D$.
If there is a set $N$ of zero capacity such that for any $z\in \partial
D\setminus N$,
\[
\lim_{D\ni x\to z}h(x)=0,
\]
then $h$ is identically zero.
\end{lemma}

\pf
Take an increasing sequence of open sets
$\{D_m\}_{m\geq 1}$ satisfying
$\overline{D_m}\subset D_{m+1}$ and
$\bigcup_{m=1}^\infty D_m =D$.
Set $\tau_m=\tau_{D_m}$. Then $\tau_m\uparrow
\tau_D$ and $\lim_{m\to\infty}X_{\tau_m}=X_{\tau_D}$ by the quasi-left
continuity of $X$.
Since $N$ has zero capacity, we have
\[
\P_x(X_{\tau_D}\in N)=0, \qquad x\in D.
\]
Therefore by the bounded convergence theorem we have for any $x\in D$,
\begin{eqnarray*}
h(x)&=&\lim_{m\to\infty}\E_x(h(X_{\tau_m}), \tau_m<\tau_D)\\
&=&\lim_{m\to\infty}\E_x(h(X_{\tau_m})1_{\partial D\setminus N}(X_{\tau_D});
\tau_m<\tau_D)=0.
\end{eqnarray*}
\qed

So far we have shown that the Martin boundary of $D$ can be
identified with a subset of the Euclidean boundary $\partial D$.

If $I$ is the set of irregular boundary points of $D$ for $X$, then $I$ is
semi-polar by Proposition II.3.3 in \cite{BG}, which is polar in our
case (Theorem 4.1.2 in \cite{FOT}). Thus $\mbox{Cap}(I)=0$. Using this observation and
the above lemma, now we can follow the proof of Theorem 4.1 in \cite{SW} and show the following
theorem, which is the main result of this section.

\begin{thm}\label{t3.1}
The Martin boundary and the minimal Martin boundary of $D$ with
respect to $X$ can be identified with the Euclidean boundary of $D$.
\end{thm}

As a consequence of Theorem \ref{t3.1}, we conclude that
for every nonnegative harmonic function $h$ for $X^D$, there exists a unique
finite measure $\mu$ on $\partial D$ such that
\begin{equation}\label{e3.1}
h(x)=\int_{\partial D}M_D(x, z)\mu(dz), \qquad x\in D.
\end{equation}
$\mu$ is called the Martin measure of $h$.

\medskip
\noindent
{\bf Acknowledgment:} The second named
author gratefully acknowledges the hospitality of Department of Mathematics of Seoul National University
where part of this work was done. The first and third named authors
gratefully acknowledge the hospitality of
Department of Mathematics of the University of Illinois at
Urbana-Champaign.

\vspace{.1in}
\begin{singlespace}

\end{singlespace}
\end{doublespace}

\end{document}